\documentclass[a4paper]{article}

\usepackage[margin=1.7in]{geometry}
\usepackage[english]{babel}
\usepackage[utf8x]{inputenc}
\usepackage{enumerate,amsmath,amssymb,amstext,amsthm}
\usepackage{graphicx}
\usepackage[colorinlistoftodos]{todonotes}
\usepackage{skull}

\usepackage{hyperref}

\newtheorem{theorem}{Theorem}[section]
\newtheorem{lemma}[theorem]{Lemma}
\newtheorem{corollary}[theorem]{Corollary}

\theoremstyle{definition}
\newtheorem{definition}[theorem]{Definition}

\newtheorem{remark}[theorem]{Remark}
\newtheorem{result}{Result}

\def\packed{\setlength{\itemsep}{1pt} \setlength{\parskip}{0pt} \setlength{\parsep}{0pt} }

\makeatletter
\renewcommand*{\@fnsymbol}[1]{\@arabic{#1}}
\makeatother

\newcommand{\R}{\mathbb{R}}

\DeclareMathOperator{\spn}{span}

\newcommand{\one}{\ensuremath{\mathbf{1}}}
\newcommand{\zero}{\ensuremath{\mathbf{0}}}

\newcommand{\chiv}{\chi_{v}}
\newcommand{\chisv}{\chi_{sv}}

\newcommand{\p}{\text{{\bf{p}}}}
\newcommand{\q}{\text{{\bf{q}}}}
\newcommand{\w}{\text{{\bf{w}}}}

\newcommand{\qeds}{\qed\vspace{.2cm}}
\DeclareMathOperator{\cone}{cone}
\DeclareMathOperator{\conv}{conv}
\DeclareMathOperator{\tr}{Tr}
\DeclareMathOperator{\sk}{sk}

\DeclareMathOperator{\rk}{rk}

\DeclareMathOperator{\im}{Im}

\title{Vector Coloring the Categorical Product of Graphs}
\author{Chris Godsil\thanks{Department of Combinatorics \& Optimization, University of Waterloo} \and David E.~Roberson\thanks{\mbox{Department of Applied Mathematics and Computer Science, Technical University of Denmark}} \and Brendan Rooney\thanks{Department of Mathematical Sciences, KAIST}
 \and Robert \v{S}\'{a}mal\thanks{Computer Science Institute, Charles University} \and Antonios Varvitsiotis\thanks{School of Physical and Mathematical Sciences, Nanyang Technological University}}

\begin{document}
\maketitle
\begin{abstract}
A vector $t$-coloring of a graph is an assignment of real vectors $p_1, \ldots, p_n$ to its vertices such that $p_i^Tp_i = t-1$ for all $i=1, \ldots, n$ and $p_i^Tp_j \le -1$ whenever $i$ and $j$ are adjacent. The vector chromatic number of $G$ is the smallest real number $t \ge 1$ for which a vector $t$-coloring of $G$ exists.  For a graph $H$ and a vector $t$-coloring $p_1,\ldots,p_n$ of a graph $G$, the assignment $(i,\ell) \mapsto p_i$ is a vector $t$-coloring of the categorical product $G \times H$. It follows that the vector chromatic number of $G \times H$ is at most the minimum of the vector chromatic numbers of the factors. We prove that equality always holds, constituting a vector coloring analog of the famous Hedetniemi Conjecture from graph coloring. Furthermore, we prove a necessary and sufficient condition for when all of the optimal vector colorings of the product can be expressed in terms of the optimal vector colorings of the factors. The vector chromatic number is closely related to the well-known Lov\'{a}sz theta function, and both of these parameters admit formulations as semidefinite programs. This connection to semidefinite programming is crucial to our work and the tools and techniques we develop could likely be of interest to others in this field.
\end{abstract}

\section{Introduction}\label{sec:intro}

For $t \ge 1$, a \emph{vector $t$-coloring} of a graph $G$ with vertex set $[n]$ is an assignment $i \mapsto p_i$ of real vectors such that 
\begin{align*}
p_i^Tp_i &= t-1 \text{ for all } i \in [n] \text{ and} \\
p_i^T p_j &\le -1 \text{ whenever } i \sim j,
\end{align*}
where `$\sim$' denotes adjacency. We note that this is different from, but equivalent to, the usual definition of a vector $t$-coloring appearing elsewhere. In the usual definition it is required that the $p_i$ are unit vectors and that $p_i^Tp_j \le -1/(t-1)$ for $i \sim j$. It is easy to see that this only differs from our definition by a scaling of the vectors (and there is no trouble with $t = 1$ in our definition). We use this slightly non-standard definition simply because it allows for a more clear connection to the semidefinite programming formulation that we will make much use of throughout the paper. Note that in a vector $1$-coloring every vertex receives the zero vector, and thus this is a valid vector coloring if and only if the graph is empty. If the graph has at least one edge, then the two vectors assigned to the ends of that edge must have inner product at most $-1$, and thus they must have norm at least $1$. This implies that the value of $t$ must be at least 2, and this can be achieved if and only if the graph is bipartite (see Section~\ref{subsec:skeletons} for details).

A vector coloring is said to be \emph{strict} if the above inequality $p_i^Tp_j \le -1$ for $i \sim j$ holds with equality for every edge. We will often use the notation $\p = (p_1, \ldots, p_n)$ to refer to the vector coloring which assigns $p_i$ to vertex $i$. The (strict) vector chromatic number of $G$, denoted $\chiv(G)$ ($\chisv(G)$), is the least real number $t$ such that $G$ admits a (strict) vector $t$-coloring. Clearly $\chiv(G) \le \chisv(G)$ by definition. Vector and strict vector colorings, as well as their associated chromatic numbers, were defined by Karger, Motwani, and Sudan in~\cite{KMS}. They noted that the strict vector chromatic number is equal to the Lov\'{a}sz $\vartheta$ number of the complement~\cite{Lovasz}, but were not aware of the fact that their vector chromatic number is equal to Schrijver's $\vartheta'$ of the complement~\cite{schrijver}. We will focus mainly on vector colorings here, but many of our results can be shown to hold for strict vector colorings without much alteration to our given proofs.

If $G$ has a $k$-coloring (in the usual sense), then mapping each color class to one of the vertices of the regular $k-1$ simplex gives a valid strict vector $k$-coloring. Thus $\chisv(G) \le \chi(G)$ and we can think of (strict) vector colorings as vector or, as we will see below, semidefinite relaxations of colorings. It is also well-known that $\omega(G) \le \chisv(G)$ where $\omega$ denotes the maximum size of a clique. This inequality remains valid for $\chiv$ in place of $\chisv$, and thus one can think of $\chiv$ as a strengthening of $\chisv$ towards $\omega$. 

Given a vector $t$-coloring $\p = (p_1, \ldots, p_n)$ of a graph $G$, we can consider the \emph{Gram matrix} of the vectors in $\p$, which we will typically denote by $M^\p$. The $ij$-entry of this matrix is equal to the inner product $p_i^Tp_j$. By the definition of vector colorings, it is easy to see that $M^\p_{ii} = t-1$ for all $i \in V(G)$ and $M^\p_{ij} \le -1$ for all $i \sim j$. Moreover, since $M^\p$ is a Gram matrix, it will be positive semidefinite. Since any positive semidefinite matrix is necessarily a Gram matrix of some set of vectors, the correspondence goes the other way as well, and so we can formulate the vector chromatic number of a graph $G$ as the following semidefinite program which already appeared (in slightly different form) in Schrijver's original paper~\cite{schrijver} defining $\vartheta'$:

\begin{equation}\label{eq:primal}
\begin{array}{lc}
\chiv(G)= & \begin{array}[t]{ll}
\min & t \tag{P}\\
\text{s.t.} & M_{ii} = t-1 \text{ for } i \in V(G) \\
 & M_{ij} \le -1 \text{ for } i \sim j \\
 & M \succeq 0
\end{array}
\end{array}
\end{equation}

Our discussion above shows that the feasible solutions to~(\ref{eq:primal}) of objective value $t$ are exactly the Gram matrices of vector $t$-colorings of $G$. Also, the dimension of the space spanned by the vectors in vector coloring $\p$, denoted $\spn(\p)$, is equal to the rank of $M^\p$. Because of this we will refer to $\dim\spn(\p)$ as the rank of $\p$. We will be interested in the maximum possible rank of an \emph{optimal} solution to (\ref{eq:primal}) for a given graph $G$, and we will refer to this maximum as the \emph{vector coloring rank of $G$}, denoted $\rk(G)$, or simply the rank of $G$ for short. If $\p$ is an optimal vector coloring of $G$ with maximum possible rank, then we will say that it is a \emph{max-rank vector coloring of $G$}, dropping the ``optimal" since we are only ever interested in the rank of optimal vector colorings. We remark that $\rk(G) = 0$ if and only if $G$ is an empty graph.

Note that considering Gram matrices of vector colorings has the effect of identifying two ``different" vector colorings $\p = (p_1, \ldots, p_n)$ and $\q = (q_1, \ldots, q_n)$ if and only if there is an isometry mapping $p_i$ to $q_i$ for all $i \in [n]$. This is analogous to identifying two classical colorings whenever they differ only by a relabelling of the colors used, i.e., they induce identical (unordered) partitions of the vertex set of the graph.

We can take the dual of the above semidefinite program in~\eqref{eq:primal} to obtain the semidefinite program in~\eqref{eq:dual} below. Here we use $i \not\simeq j$ to denote that $i$ and $j$ are neither equal nor adjacent, and $\text{sum}(B)$ to denote the sum of the entries of the matrix $B$ (this is also equal to $\tr(BJ)$ where $J$ is the all-ones matrix). We note that this dual also originally appeared in Schrijver's paper~\cite{schrijver}.

\begin{equation}\label{eq:dual}
\begin{array}{lc}
\chiv(G) \ = & \begin{array}[t]{ll}
\max & \text{sum}(B) \tag{D} \\
\text{s.t.} & B_{ij} = 0 \text{ for } i \not\simeq j \\
 & B_{ij} \ge 0 \text{ for all } i,j \\
 & \tr(B) = 1 \\
 & B \succeq 0
\end{array}
\end{array}
\end{equation}

Note that both of these semidefinite programs are bounded and strictly feasible, thus Slater's condition holds and so they have the same optimal values. We will refer to feasible (optimal) solutions to~(\ref{eq:primal}) and~\eqref{eq:dual} as feasible (optimal) \emph{primal} and \emph{dual} solutions for $\chiv(G)$ respectively. \\

Our primary focus in this work is on vector colorings of categorical products of graphs. Given graphs $G$ and $H$, their \emph{categorical product}, denoted $G \times H$, has vertex set $V(G) \times V(H)$ where vertex $(i,\ell)$ is adjacent to $(j,k)$ if $i \sim j$ and $\ell \sim k$. Given graphs $G$ and $H$, and a $c$-coloring $\varphi$ of $G$, one can construct a $c$-coloring of the categorical product $G \times H$ by mapping $(i,\ell)$ to $\varphi(i)$ for all $i \in V(G)$, $\ell \in V(H)$. We say that such a $c$-coloring of $G \times H$ is \emph{induced by $\varphi$}, or simply \emph{induced by $G$}. Note that a coloring of $G \times H$ is induced by $G$ if and only if the color of a vertex does not depend on its $H$-coordinate (unless $H$ is an empty graph and $G$ is not). Since both $G$ and $H$ can induce colorings of $G \times H$, we see that
\begin{equation}\label{eq:hedineq}
\chi(G \times H) \le \min\{\chi(G), \chi(H)\}.
\end{equation}
In 1966, Hedetniemi conjectured that equality holds, and only a few special cases have been proven in the time since then. Most notably, El-Zahar and Sauer~\cite{hedetniemi4} gave a proof for when the minimum is four (for smaller values of the minimum, the proof is straightforward), but for all larger values the conjecture remains open. More recently, it was shown by Zhu~\cite{frached} that the conjecture holds when chromatic number is replaced by fractional chromatic number.\\

As in the case of colorings, vector colorings of $G$ and $H$ induce vector colorings of the categorical product $G \times H$. Concretely, given graphs $G$ and $H$ and a vector $t$-coloring $\p = (p_1, \ldots, p_n)$ of $G$, the map $(i,\ell) \mapsto p_i$ is easily seen to be a vector $t$-coloring of $G \times H$. We say that this vector coloring is induced by $\p$, or simply that it is induced by $G$. As with colorings, a vector coloring of $G \times H$ is induced by $G$ if and only if\footnote{As in the coloring case, there is an exception when $H$ is empty since then $G \times H$ is empty and the only optimal vector coloring assigns the zero vector everywhere. In this case the vector assigned to a vertex does not depend on its $H$-coordinate, but if $G$ is not empty then this vector coloring cannot be realized as one induced by $G$.} the vector assigned to $(i,\ell)$ does not depend on $\ell$. If $M^\p$ is the Gram matrix of the vector coloring $\p$, then the Gram matrix of the vector coloring of $G \times H$ induced by $\p$ is given by the Kronecker product $M^\p \otimes J$. Equivalently, a vector coloring of $G \times H$ is induced by $G$ if and only if the $(i,\ell)(j,k)$-entry of its Gram matrix only depends on $i$ and $j$.

This work is motivated by the following two questions:

\begin{enumerate}
\item When are the vector colorings induced by the factors optimal for $G \times H$?
\item When is it possible to describe all optimal vector colorings of $G \times H$ in terms of the vector colorings of the individual factors?
\end{enumerate}

\subsection{Summary of results}

Since both $G$ and $H$ can induce vector colorings of $G \times H$, it follows that $\chiv(G \times H) \le \min\{\chiv(G), \chiv(H)\}$. Our first result in this work is the vector coloring analog of Hedetniemi's Conjecture, i.e.,

\begin{result}\label{res:veched}
For any graphs $G$ and $H$, we have that
\[\chiv(G \times H) = \min\{\chiv(G), \chiv(H)\}.\]
\end{result}
The analogous result for strict vector colorings was recently proven by Severini and a subset of the authors in~\cite{sabvshed}. However, the technique used there does not extend to the vector chromatic number. On the other hand, our proof of the above can be adapted to give a shorter proof of the result for the strict vector chromatic number (see Theorem~\ref{thm:hedetniemis}).

Result~\ref{res:veched} answers the first question from above, showing that the optimal vector colorings induced by the factor(s) with minimum vector chromatic number always induce optimal vector colorings of the product. Motivated by this, we aim to describe \emph{all} of the optimal vector colorings of $G \times H$ in terms of vector colorings of the factors, i.e., to answer the second question above. We distinguish two cases.

First, in the case where $\chiv(G) < \chiv(H)$, only the vector colorings induced by $G$ are optimal, not those induced by $H$. Further, note that the rank of a vector coloring induced by $G$ is merely the rank of the corresponding vector coloring of $G$. Thus, in order for every optimal vector coloring of $G \times H$ to be induced by $G$, a trivial necessary condition is that $\rk(G \times H) = \rk(G)$. Surprisingly, our second main result shows that this sufficient as well.

\begin{result}\label{res:B}
Let $G$ and $H$ be graphs with $\chiv(G) < \chiv(H)$. Then every optimal vector coloring of $G \times H$ is induced by $G$ if and only if $\rk(G \times H) = \rk(G)$.
\end{result}

Second, in the case where both factors have the same vector chromatic number, i.e., $\chiv(G) = \chiv(H)$, each factor can induce optimal vector colorings of the product. The Gram matrix of such an induced vector coloring has the form $M^\p \otimes J$ or $J \otimes M^\q$, where $M^\p$ and $M^\q$ are Gram matrices of optimal vector colorings $\p$ and $\q$ of $G$ and $H$ respectively. However, one can also form convex combinations of the Gram matrices of vector colorings induced by each of the factors. This results in an optimal vector coloring of $G \times H$ whose Gram matrix has the form $\alpha (M^\p \otimes J) + \beta (J \otimes M^\q)$ for some $0 \le \alpha = 1-\beta \le 1$. We will abuse terminology somewhat and refer to any vector coloring whose Gram matrix has this form as a \emph{convex combination} of the vector colorings induced by $\p$ and $\q$. Note that a particular vector coloring with this Gram matrix is given by $(i,\ell) \mapsto \sqrt{\alpha} p_i \oplus \sqrt{\beta} q_\ell$ for all $i \in V(G)$, $\ell \in V(H)$. We refer to this as a \emph{direct sum} of $\p$ and $\q$ and denote it by $\sqrt{\alpha} \p \oplus \sqrt{\beta} \q$. Such a mixing of vector colorings induced by the factors is more difficult to recognize, and indeed we do not know a simple necessary and sufficient condition for when a vector coloring has this form. This added difficulty, along with other considerations, causes an increase in the complexity of our arguments for this case.

Note that since $J$ has rank one, the rank of $\alpha (M^\p \otimes J) + \beta (J \otimes M^\q)$ is at most $\rk(M^\p) + \rk(M^\q)$. Thus the maximum possible rank of an optimal vector coloring of $G \times H$ which is a convex combination of vector colorings induced by the factors is $\rk(G) + \rk(H)$. It is not immediate that this rank can be realized by such a vector coloring, but it only requires a short proof. Analogously to the previous case, this implies that if every optimal vector coloring of $G \times H$ is a convex combination of vector colorings induced by the factors, then we must have that $\rk(G \times H) \le \rk(G) + \rk(H)$. Also analogous to the previous case, this turns out to be sufficient:

\begin{result}\label{res:A}
Let $G$ and $H$ be graphs with $\chiv(G) = \chiv(H)$. Then every optimal vector coloring of $G \times H$ is a convex combination of vector colorings induced by $G$ and $H$ if and only if $\rk(G \times H) = \rk(G) + \rk(H)$.
\end{result}

We remark that it is not always the case that every optimal vector coloring of $G \times H$ is a convex combination of vector colorings induced by $G$ and $H$ when $\chiv(G) = \chiv(H)$ (and analogously in the $\chiv(G) < \chiv(H)$ case), and so the assumptions on $\rk(G \times H)$ in Results~\ref{res:B} and~\ref{res:A} are not superfluous.

\subsection{Applications}\label{subsec:applications}

In order to apply Results~\ref{res:B} and~\ref{res:A}, we need to determine the maximum rank of an optimal vector coloring for each of $G$, $H$, and $G \times H$. This task is difficult in general, but we address this by using the duality theory of semidefinite programming.

Let $M$ and $B$ be primal/dual feasible solutions for $\chiv(G)$. Semidefinite programming duality theory implies that if $(M,B)$ are primal/dual optimal, then $MB = 0$ (for details see Lemma~\ref{lem:compslack}). Furthermore, this shows that if $(M,B)$ are primal/dual optimal solutions, then $\rk(M) \le \text{corank}(B)$, where $\text{corank}(B)$ denotes the dimension of the kernel/null space of $B$. Thus any optimal dual solution $B$ provides an upper bound on the maximum rank of an optimal vector coloring. We say that a pair of primal/dual optimal solutions $(M,B)$ are \emph{strictly complementary} if $\rk(M) = \text{corank}(B)$. Note that in this case, $\rk(M) = \rk(G)$, so $B$ serves as a certificate that $M$ has the largest possible rank. Lastly, we say that $G$ satisfies strict complementarity if there exists a pair of primal/dual optimal strictly complementary solutions.

Roughly speaking, to prove that $\rk(G \times H) \le \rk(G)$ (in the case where $\chiv(G) < \chiv(H)$) we show that if $G$ satisfies strict complementarity, then $G \times H$ also satisfies strict complementarity. Similarly, to prove $\rk(G \times H) \le \rk(G) + \rk(H)$ (in the case where $\chiv(G) = \chiv(H)$) we show that if both $G$ and $H$ satisfy strict complementarity, then so does $G \times H$. This is a remarkable property of the vector chromatic number, and an interesting research direction is to find other classes of semidefinite programs that enjoy this property.

Using Result~\ref{res:B} and the preceding discussion, in Corollary~\ref{cor:B} we prove:

\begin{result}\label{res:Bcor}
Let $G$ and $H$ be graphs such that $\chiv(G) < \chiv(H)$ and $H$ is connected. If $G$ admits a strictly complementary dual solution with strictly positive diagonal, then every optimal vector coloring of $G \times H$ is induced by $G$.
\end{result}

This result has some interesting consequences. First, it allows us to find classes of graphs for which we can describe all optimal vector colorings of their categorical products. Specifically, a graph $G$ is 1-walk-regular if for all $k \in \mathbb{N}$,
\begin{enumerate}
\item the number of walks of length $k$ starting and ending at a vertex of $G$ is independent of the vertex;
\item the number of walks of length $k$ starting at one end of an edge and ending at the other is independent of the edge.
\end{enumerate}

Note that any 1-walk-regular graph must be regular. Also, any graph which is vertex- and edge-transitive is easily seen to be 1-walk-regular. Other classes of 1-walk-regular graphs include distance regular graphs and, more generally, graphs which are a single class in an association scheme.

It follows from results in~\cite{UVC1} that any 1-walk-regular graph has a strictly complementary dual solution with strictly positive diagonal. As a consequence, if $G$ is 1-walk-regular and $H$ is connected with $\chiv(G) < \chiv(H)$, then all of the optimal vector colorings of $G \times H$ are induced by $G$.

Result~\ref{res:Bcor} also generalizes a result of Pak and Vilenchik~\cite{Pak}. Specifically, they showed that if an $r$-regular graph $H$ with eigenvalues $\lambda_1 \ge \lambda_2 \ge \ldots \ge \lambda_n$ satisfies $\lambda(H) < r/(m-1)$, where $\lambda(H) = \max_{i \ge 2} |\lambda_i |$, then the product $K_m \times H$ has a unique vector $m$-coloring (the one induced by $K_m$). It turns out that their assumptions on $H$ imply that it is both connected and $\chiv(H) > m = \chiv(K_m)$, and so their result is a special case of Result~\ref{res:Bcor} (it is not difficult to show that $K_m$ is satisfies strict complementarity and has a unique vector $m$-coloring). The details of this are given in Section~\ref{subsec:Btheorem}. We remark that the result of Pak and Vilenchik was the original inspiration for the research presented in this work, particularly in the case of Results~\ref{res:B} and~\ref{res:Bcor}.

Analogously, based on Result~\ref{res:A}, in Corollary~\ref{cor:A} we prove the following, where we say that a matrix is \emph{connected} if its underlying graph is connected (see Section~\ref{subsec:prelims}):

\begin{result}\label{res:Acor}
Let $G$ and $H$ be graphs such that $\chiv(G) = \chiv(H)$. If both $G$ and $H$ admit connected strictly complementary dual solutions, then every optimal vector coloring of $G \times H$ is a convex combination of vector colorings induced by the factors.
\end{result}

\subsection{Motivations from graph coloring}\label{subsec:motivations}

The other main motivation for our investigations was the work of Duffus, Sands, and Woodrow. In~\cite{duffus}, they uncovered a connection between unique colorability and Hedetniemi's Conjecture. A graph $G$ is \emph{uniquely $c$-colorable} if $V(G)$ has a unique partition into at  most $c$ nonempty independent sets, i.e., $\chi(G) = c$ and $G$ has a unique $c$-coloring up to relabeling of the colors. Duffus, Sands, and Woodrow considered the following three parameterized statements:
\begin{itemize}
\item $\boldsymbol{(A_n):}$ For all uniquely $n$-colorable graphs $G$ and $H$, each $n$-coloring of $G \times H$ is induced by $G$ or $H$.
\item $\boldsymbol{(B_n):}$ For all uniquely $n$-colorable graphs $G$ and connected graphs $H$ with $\chi(H) > n$, the graph $G \times H$ is uniquely $n$-colorable.
\item $\boldsymbol{(C_n):}$ For all graphs $G$ and $H$ with $\chi(G) = \chi(H) = n$, we have
\[\chi(G \times H) = n.\]
\end{itemize}
Note that, since any graph with chromatic number at least $n$ contains a subgraph with chromatic number exactly $n$, Hedetniemi's Conjecture is equivalent to $(C_n)$ being true for all $n \in \mathbb{N}$. Surprisingly, Duffus, Sands, and Woodrow showed that $(A_n) \Rightarrow (B_n) \Rightarrow (C_{n+1})$ for all $n$. Unfortunately, they were not able to prove $(A_n)$ or $(B_n)$ in general, but only under additional restrictions, such as one of the factors being a complete graph.

It is possible to define unique vector colorability as well, which was originally done in~\cite{Pak}. Given any vector $t$-coloring $\p$ of a graph $G$, it is easy to see that applying any orthogonal transformation to the vectors in $\p$ will produce a vector $t$-coloring of $G$ (which will essentially always be different from $\p$). So in this sense, there is never a unique vector $t$-coloring of a graph (unless it is empty). Of course, this is analogous to relabelling the colors in a classical coloring, and so we merely need to quotient out by this equivalence. As mentioned earlier, we do this by considering the Gram matrices of vector colorings. Thus we have the following definition:

\begin{definition}
A graph $G$ is \emph{uniquely vector colorable} if for any two \emph{optimal} vector colorings $\p = (p_1, \ldots, p_n)$ and $\q = (q_1, \ldots, q_n)$, we have that
\[p_i^Tp_j = q_i^Tq_j \text{ for all } i,j \in V(G).\]
Equivalently, there is a unique optimal primal solution for $\chiv(G)$. We say that a graph is \emph{uniquely vector $t$-colorable} if it is uniquely vector colorable with vector chromatic number $t$.
\end{definition}

We can now devise vector coloring analogs of the three statements considered by Duffus, Sands, and Woodrow:

\begin{itemize}
\item $\boldsymbol{(A'):}$ For all uniquely vector colorable graphs $G$ and $H$ with $\chiv(G) = \chiv(H)$, every optimal vector coloring of $G \times H$ is a convex combination of the vector colorings induced by $G$ and $H$.
\item $\boldsymbol{(B'):}$ For all uniquely vector colorable graphs $G$ and connected graphs $H$ with $\chiv(G) < \chiv(H)$, the graph $G \times H$ is uniquely vector colorable.
\item $\boldsymbol{(C'):}$ For all graphs $G$ and $H$,
\[\chiv(G \times H) = \min\{\chiv(G), \chiv(H)\}.\]
\end{itemize}

It is not difficult to see the similarities between the statements $(A')$, $(B')$, and $(C')$ and our Results~\ref{res:A}, \ref{res:B}, and~\ref{res:veched} (in that order). Indeed, Result~\ref{res:veched} is exactly statement $(C')$, whereas $(A')$ and $(B')$ are Results~\ref{res:A} and~\ref{res:B} in the case of uniquely vector colorable $G$ (and $H$ in the former), but without the assumptions on $\rk(G \times H)$. Thus the statements $(A')$ and $(B')$ are true if and only if the assumptions on $\rk(G \times H)$ from Results~\ref{res:A} and~\ref{res:B} must hold whenever $G$ (and $H$ in the former) are uniquely vector colorable. Unfortunately, we do not yet know how to show this, so it remains an interesting open problem.

\subsection{Outline}

The rest of the paper is outlined as follows. In Section~\ref{sec:SDPs}, we prove that every optimal vector coloring of a graph can be obtained in a specific manner from a max-rank vector coloring. This is one of the main tools we use for our results in this work. In Section~\ref{subsec:duality}, we prove the complementary slackness conditions for the primal/dual pair of semidefinite programs given in~\eqref{eq:primal} and~\eqref{eq:dual}, and we review the notion of strictly complementary pairs of solutions. Following this, in Section~\ref{subsec:certification} we introduce another formulation for vector chromatic number for which it is easy to combine solutions for two graphs to construct a solution for their product. Next, in Section~\ref{subsec:hedetniemi}, we use this reformulation to prove the vector coloring analog of Hedetniemi's Conjecture, i.e., Result~\ref{res:veched}. Section~\ref{sec:techlems} introduces and develops the concepts of skeletons and neighborliness in vector colorings. We prove several lemmas about these notions that are crucial to the proof of Result~\ref{res:A}. Section~\ref{sec:products} contains the proofs of Results~\ref{res:B}--\ref{res:Acor}.We also provide some conditions on the skeletons of the factors which are necessary for the only optimal vector colorings of the product to be the convex combinations of the vector colorings induced by the factor(s) with the minimum vector chromatic number. In Section~\ref{subsec:1walkreg}, we show that statements $(A')$ and $(B')$ hold for (connected) 1-walk-regular graphs even without the assumption of unique vector colorability. In Section~\ref{subsec:implications}, we prove a vector coloring analog of the result of Duffus, Sands, and Woodrow that $(A_n) \Rightarrow (B_n)$ for all $n \in \mathbb{N}$. Finally, in Section~\ref{sec:discuss} we discuss our results and some possible open questions.

\subsection{Preliminaries and basic notation}\label{subsec:prelims}

Throughout we set $[n]=\{1,\ldots,n\}$. We denote by $\one$ the all-ones vector and by $\zero$ the all-zeros vector of appropriate size. All vectors are column vectors. Furthermore, we denote by $\spn(p_1,\ldots,p_n)$ the linear span of the vectors $\{p_i\}_{i=1}^n$.   
 
A {\em convex cone} in $\mathbb{R}^d$ is a subset of vectors that is closed under positive linear combinations. Given a finite set of vectors $S = \{p_1, \ldots, p_k\} \subseteq \mathbb{R}^d$, the \emph{conical hull} of $S$, denoted $\cone(S)$, is given by
\[\cone(S) = \left\{\sum_{i = 1}^k \alpha_i p_i : \alpha_i \ge 0\right\},\]
and is 
 always a closed convex cone. The \emph{convex hull} of a finite set of vectors $S = \{p_1, \ldots, p_k\} \subseteq \mathbb{R}^d$, denoted $\conv(S)$, is the set of all convex combinations of the vectors in $S$, i.e.,
\[\conv(S) = \left\{\sum_{i=1}^k \alpha_i p_i : \sum_{i=1}^k \alpha_i = 1 \ \& \ \alpha_i \ge 0 \text{ for all } i\in [k]
\right\}.\]

Given a $n\times n$ matrix $X$ we denote its kernel/null space by $\ker X$  and its image/column space by $\im X$. We will use $\lambda_{\max}(X)$ and $\lambda_{\min}(X)$ to denote the maximum and minimum eigenvalues of $X$ respectively. A symmetric matrix is {\em positive semidefinite} if all of its eigenvalues are nonnnegative. The {\em Gram matrix} of a set of vectors $p_1, \ldots, p_n$ is the $n\times n$ matrix with $ij$-entry equal to $p_i^T p_j$. This matrix is positive semidefinite and its rank is equal to $\dim \spn(p_1,\ldots,p_n)$. Conversely, a $n\times n$ positive semidefinite matrix with rank equal to $r$ can always be realized as the Gram matrix of a family of real vectors $p_1,\ldots,p_n\in \R^r$. 

The {\em Schur} product of two $n\times n$ matrices $X,Y$, denoted by $X\circ Y$, is the $n\times n$ matrix whose entries are given by $(X\circ Y)_{ij}=X_{ij}Y_{ij}$ for all $i,j\in [n]$. The Schur product of two positive semidefinite matrices is positive semidefinite. Let $X, Y$ be matrices with dimensions $a \times b$ and $c \times d$ respectively. The {\em Kronecker} product of $X$ and $Y$, denoted by $X\otimes Y$, is the  $ac\times bd$ block matrix
\[\begin{pmatrix}X_{11}Y & \ldots & X_{1b}Y\\
\vdots &  & \vdots\\
X_{a1}Y& \ldots & X_{ab}Y\end{pmatrix}.\]
The Kronecker product of two positive semidefinite matrices is also positive semidefinite. Furthermore, if $\lambda_1,\ldots,\lambda_n$  and $\mu_1,\ldots,\mu_m$ are the eigenvalues of $X$ and $Y$ respectively, the eigenvalues of $X\otimes Y$ are $\lambda_i\mu_j, i\in [n], j\in [m].$ Lastly, the sum of two positive semidefinite matrices $X,Y$ is also positive semidefinite and furthermore ${\rm rank}(X+Y)\ge {\rm max}\{{\rm rank}(X),{\rm rank}(Y)\}.$

The {\em graph} of a symmetric $n\times n$ matrix $M$, denoted by $G(M)$, is the graph with vertex set $[n]$ where distinct vertices $i$ and $j$ are adjacent if $M_{ij}\ne 0$. We say that a  matrix is \emph{connected} if the graph $G(M)$ is connected in the graph-theoretic sense. If $M$ is positive semidefinite, a necessary condition for $G(M)$ to be connected is that all diagonal entries are non-zero. 

The Perron-Frobenius Theorem states that the maximum eigenvalue of a connected, entrywise nonnegative matrix is also maximum in absolute value, has multiplicity 1, and admits an entrywise positive eigenvector. If a nonnegative matrix $M$ is not connected, then it is, up to a common permutation of its rows and columns, a direct sum of connected nonnegative matrices (these correspond to the components of $G(M)$). It follows that any nonnegative matrix has an entrywise nonnegative eigenvector for its maximum eigenvalue.

\section{Characterizing Optimal Vector Colorings}\label{sec:SDPs}

In this section we show that we can describe all optimal vector colorings of a graph $G$ provided that we know one max-rank optimal vector coloring. In order to prove this we first need the following two lemmas, the first of which was proven in~\cite{UVC1} and so we omit the proof.

\begin{lemma}\label{lem:R}
Let $P$ be an $n \times d$ matrix with rank $d$. If $\im (Y) \subseteq \im(PP^T)$, then there exists a symmetric $d \times d$ matrix $R$ such that $Y=PRP^T$.
\end{lemma}

The proof  of Lemma \ref{lem:R} can be found in~\cite{UVC1}. Next we show that  max-rank optimal solutions for $\chiv(G)$ have the largest image among all optimal solutions. This a well-known property of semidefinite programs, e.g.~see \cite[Lemma 2.3]{deklerk}. We give a proof for completeness.

\begin{lemma}\label{lem:image}
Let $G$ be a graph and let $\hat{M}$ be an optimal primal solution for $\chiv(G)$ with maximum possible rank. Then $\im(\hat{M}) \supseteq \im(M)$ for any optimal primal solution $M$ for $\chiv(G)$.
\end{lemma}
\proof
Let $\hat{M}$ be as in the lemma statement and suppose that the conclusion does not hold for some optimal solution $M$. Define $M' = (1/2)(\hat{M} + M)$. It is obvious that $M'$ is an optimal primal solution for $\chiv(G)$. We will show that $M'$ has strictly greater rank than $\hat{M}$.

Since $M'$ is positive semidefinite, a vector $p$ is in its kernel if and only if $p^TM'p = 0$. But since each of $\hat{M}$ and $M$ are positive semidefinite, this happens if and only if $p^T\hat{M}p = 0$ and $p^TMp = 0$. In other words, $\ker(M') = \ker(\hat{M}) \cap \ker(M)$. By assumption, $\im(\hat{M}) \not\supseteq \im(M)$ and therefore $\ker(\hat{M}) \not\subseteq \ker(M)$. Thus $\ker(M')$ is strictly smaller than $\ker(\hat{M})$ and it follows that $M'$ has strictly larger rank than $\hat{M}$.\qeds

Now we can prove the main result of this section:

\begin{theorem}\label{thm:main}
Suppose that $\p = (p_1, \ldots, p_n)$ is a max-rank vector coloring of $G$ with $p_i \in \mathbb{R}^{\rk(G)}$. Let $P$ be the matrix whose $i^\text{th}$ row is $p_i^T$. Then $M$ is the Gram matrix of an optimal vector coloring of $G$ if and only if

\[M = P(I+R)P^T,\]
for some $\rk(G) \times \rk(G)$ symmetric matrix $R$ satisfying
\begin{enumerate}[(i)]\packed
    \item $p_i^T R p_i = 0 \text{ for all } i \in [n]$;
    \item $p_i^T R p_j \le -1 - p_i^Tp_j \text{ for } i \sim j$;
    \item $I+R \succeq 0$.
\end{enumerate}
\end{theorem}
\proof
Suppose that $M = P(I+R)P^T$ and $R$ satisfies the conditions stated in the theorem. Since $I+R \succeq 0$, we have that $M \succeq 0$ and thus is the Gram matrix of some set of vectors $q_i$ for $i \in V(G)$. Letting $t = \chiv(G)$, we also have that $M_{ii} = (PP^T)_{ii} + (PRP^T)_{ii} = p_i^Tp_i + p_i^TRp_i = t-1 + 0 = t-1$, by the first condition on $R$. Lastly, for $i \sim j$,
\[M_{ij} = p_i^Tp_j + p_i^T R p_j \le p_i^Tp_j + \left(-1 - p_i^Tp_j\right) = -1.\]
Therefore, $M$ is the Gram matrix of an optimal vector coloring of $G$.

Conversely, suppose that $M$ is the Gram matrix of an optimal vector coloring of $G$. By Lemma~\ref{lem:image}, we have that $\im(M) \subseteq \im (PP^T)$. Let $Y = M - PP^T$ and note that $\im(Y) \subseteq \spn(\im(M) \cup \im(PP^T)) = \im(PP^T)$. Thus by Lemma~\ref{lem:R} there exists a symmetric matrix $R$ such that $Y = PRP^T$. Therefore, $M = Y + PP^T = PRP^T + PP^T = P(I+R)P^T$.

Note that $p_i^TRp_j = (PRP^T)_{ij} = M_{ij} - (PP^T)_{ij}$ for all $i,j$. Since both $PP^T$ and $M$ are Gram matrices of optimal vector colorings of $G$ by assumption, we have that $(PP^T)_{ii} = M_{ii} = t-1$ for all $i$. Therefore, we have that
\[p_i^TRp_i = 0,\]
as required. Furthermore, for $i \sim j$ we have that
\[p_i^TRp_j = M_{ij} - (PP^T)_{ij} \le -1 - p_i^Tp_j.\]
Lastly, since $P$ has full column rank, we have that $M = P(I+R)P^T$ is positive semidefinite if and only if $I+R$ is positive semidefinite. Since $M$ is a Gram matrix, we have that $I+R \succeq 0$ as required.\qeds

Note that the max-rank condition on $\p$ in the above is necessary in the following sense: the matrix $P(I+R)P^T$ has rank at most that of $P$ regardless of $R$, thus we can only obtain Gram matrices of vector colorings of equal or lesser rank using the construction from the theorem. So if we start with a vector coloring that is not of maximum rank, then we can never obtain a max-rank vector coloring from this construction.

\begin{remark}\label{rem:earlymain}
Theorem 2.3 from our earlier work~\cite{UVC1} is similar to Theorem~\ref{thm:main} here, but the former concerns more general assignments of vectors to the vertices of graphs. However, we note that in the case of vector colorings the above theorem is actually stronger. This is because the theorem in~\cite{UVC1} required the existence of a dual certificate (called a ``spherical stress matrix"), and furthermore required that the max-rank vector coloring $\p$ be an optimal \emph{strict} vector coloring as well. Since then we have discovered that these additional assumptions are superfluous and that the essential property of $\p$ is that it has maximum possible rank.
\end{remark}

Note that in practice, when searching for a matrix $R$ satisfying the conditions in Theorem~\ref{thm:main} above, it suffices to find $R$ such that $p_i^TRp_j = 0$ whenever $i = j$ and whenever $i \sim j$ and $p_i^Tp_j = -1$. Given such an $R$ it is always possible to scale it so that all of the other inequalities of $(ii)$ hold and such that $I+R \succeq 0$. Of course, if one can show that no such matrix $R$ exists, then this proves that the graph is uniquely vector colorable. This is the technique used to prove unique vector colorability in our previous works~\cite{UVC1} and~\cite{UVC2}. These works also present efficient algorithms for finding matrices $R$ satisfying the conditions of the above theorem.

\subsection{Duality and complementary slackness}\label{subsec:duality}

In order to apply Theorem~\ref{thm:main} we must first obtain a max-rank vector coloring of our graph. The difficulty here is not so much in finding an optimal vector coloring, but rather verifying that it is of max-rank. We do not know of any general method for doing this, but here we will present a sufficient condition for an optimal vector coloring to be of maximum rank.

The main tool we will use is complementary slackness for semidefinite programs. For the primal/dual pair of semidefinite programs in~\eqref{eq:primal} and~\eqref{eq:dual}, the complementary slackness conditions are given in the lemma below. We remark that this is a standard result in the theory of semidefinite programs, but we give a proof for completeness.

\begin{lemma}\label{lem:compslack}
Let $G$ be a graph and let $M$ and $B$ be feasible primal and dual solutions for $\chiv(G)$ respectively. If $M$ and $B$ have objective values $t$ and $s$, then
\[\tr(MB) = (t - s) + \sum_{i \sim j} (M_{ij}+1)B_{ij}.\]
In particular, this implies that $M$ and $B$ are both optimal if and only if $MB = 0$ and $M_{ij} < -1 \Rightarrow B_{ij} = 0$ for all $i,j$.
\end{lemma}
\proof
Defining $\tilde{M} = M - tI + J$, we have that $M = tI - J + \tilde{M}$. Therefore,
\begin{align*}
\tr(MB) &= \tr((tI - J + \tilde{M})B) \\
&= t \tr(B) - \tr(BJ) + \tr(\tilde{M}B) \\
&= (t - s) + \sum_{ij} \tilde{M}_{ij}B_{ij} \\
&= (t-s) + \sum_{i \sim j} \tilde{M}_{ij}B_{ij} \\
&= (t-s) + \sum_{i \sim j} (M_{ij} + 1)B_{ij}.
\end{align*}
This proves the first claim of the lemma. Now note that the summation in the last expression is always non-positive, since $M_{ij} \le -1$ for all $i \sim j$. Moreover, this summation is zero if and only if $M_{ij} < -1 \Rightarrow B_{ij} = 0$ for all $i,j$. Now suppose that $M$ and $B$ are both optimal. Then by strong duality we have that $t = s$ and therefore $\tr(MB) = \sum_{i \sim j} (M_{ij} + 1)B_{ij} \le 0$. However, since both $M$ and $B$ are positive semidefinite, we have that $\tr(MB) \ge 0$. Therefore, we have that $\tr(MB) = \sum_{i \sim j} (M_{ij} + 1)B_{ij} = 0$. This further implies that $MB = 0$ and $M_{ij} < -1 \Rightarrow B_{ij} = 0$ for all $i,j$.

Conversely, if $MB = 0$ and $M_{ij} < -1 \Rightarrow B_{ij} = 0$ for all $i,j$, then
\[0 = \tr(MB) = (t-s) + \sum_{i \sim j} (M_{ij} + 1)B_{ij} = t-s,\]
and thus $t = s$, i.e., both $M$ and $B$ are optimal.\qeds

\begin{remark}\label{rem:G(B)}
It will be useful to consider the contrapositive of the condition $M_{ij} < -1 \Rightarrow B_{ij} = 0$ for all $i,j$. This is of course $B_{ij} \ne 0 \Rightarrow M_{ij} \ge -1$ for all $i,j$. By the feasibility conditions on $M$ and $B$, this is equivalent to $B_{ij} > 0 \Rightarrow M_{ij} = -1$ for all $i \sim j$.
\end{remark}

As mentioned in Section~\ref{subsec:applications}, a consequence of the above complementary slackness conditions is that for any pair of primal/dual optimal solutions $M$ and $B$, we have that $\im(M) \subseteq \ker(B)$ and analogously $\ker(M) \supseteq \im(B)$. This implies that $\rk(M) + \rk(B) \le |V(G)|$, and we will say that $M$ and $B$ are \emph{strictly complementary} if equality holds, which is equivalent to $\im(M) = \ker(B)$. We will also say that $G$ \emph{satisfies strict complementarity (for vector chromatic number)} if such a pair of strictly complementary solutions exists. Note that since $\rk(M) + \rk(B) \le |V(G)|$ holds for any pair of optimal solutions, if there exists some pair of strictly complementary solutions then a particular pair of optimal solutions $M'$ and $B'$ are strictly complementary if and only if both $M'$ and $B'$ have the maximum possible rank. In this case $\text{corank}(B') = \rk(M') = \rk(G)$. Thus we have the following lemma:

\begin{lemma}\label{lem:mrcert}
Let $G$ be a graph with optimal vector coloring $\p$. Then $\p$ has maximum rank if there exists an optimal dual solution $B$ such that $\text{corank}(B) = \dim\spn(\p)$.
\end{lemma}

Since an optimal dual solution $B$ is part of a strictly complementary pair of solutions if and only if $\text{corank}(B) = \rk(G)$, we will sometimes refer to a dual solution $B$ with this property as a \emph{strictly complementary dual solution} without explicitly mentioning a corresponding primal solution.

\subsection{A useful reformulation for $\chiv$}\label{subsec:certification}

We have seen that vector colorings of the factors of a categorical product induce vector colorings of the product. In other words, we can use primal solutions for $\chiv(G)$ and $\chiv(H)$ to build primal solutions for $\chiv(G \times H)$. However, we will also need a way to use dual solutions for the factors to find dual solutions for the product. A first approach may be to take Kronecker products of dual solutions for the factors in order to obtain a dual solution for the product. But it is not hard to see that this does not work. Indeed, the objective value of such a solution would be the product of the objective values, whereas we know that $\chiv(G \times H)$ is at most the minimum of $\chiv(G)$ and $\chiv(H)$.

It turns out that the right approach is to use another formulation for vector chromatic number which works well with the categorical product. This formulation appeared in~\cite{galtman}, and an analogous formulation for Lov\'{a}sz theta appeared even in Lov\'{a}sz' original paper~\cite{Lovasz}. We present the formulation in the lemma below along with the standard proof, since converting between feasible solutions for it and feasible solutions for the dual~(\ref{eq:dual}) will be important for our later results. 
Note that below, $i \not\sim j$ includes the case $i = j$. Also, we use $\|M\|$ to denote the maximum (absolute value) of the eigenvalues of the matrix $M$.

\begin{lemma}\label{lem:dualconvert}
For any graph $G$,
\begin{equation}\label{eq:Adual}
\begin{array}{lc}
\chiv(G)= & \begin{array}[t]{ll}\tag{D$'$}
\max & \|I + A\| \\
\text{s.t.} & A_{ij} = 0 \text{ if } i \not\sim j \\
 & A_{ij} \ge 0 \text{ for all } i,j\\
 & I + A \succeq 0
\end{array}
\end{array}
\end{equation}
\end{lemma}
\proof
Suppose that $B$ is a feasible dual solution for $\chiv(G)$. For ease of presentation, we will assume that $B_{ii} > 0$ for all $i \in V(G)$, but the proof is easy to adapt to the general case. Let $D$ be the diagonal part of $B$, i.e.,~$D = I \circ B$. Note that by our assumption on $B$, the matrix $D$ has only positive diagonal entries and is thus invertible and has a square root. We will show that $A =  D^{-1/2}BD^{-1/2} - I$ is a feasible solution to~(\ref{eq:Adual}) for $G$ with objective value at least that of $B$. First note that since $B \succeq 0$, we have that $I+A = D^{-1/2}BD^{-1/2} \succeq 0$. Similarly, since the diagonal entries of $D$ are positive, we have that all entries of $I+A$ are nonnegative. Furthermore, it is easy to see that $B$ and $D^{-1/2}BD^{-1/2}$ have the same zero entries, and the latter has 1's on the diagonal. Therefore, $A$ has 0's in all the required positions to be a feasible solution to~(\ref{eq:Adual}). To see that $A$ has objective value at least $\text{sum}(B)$, let $v$ be a vector with $v_i = (B_{ii})^{1/2}$. Then, since $\tr(B) = 1$, we have that $v$ is a unit vector, and moreover $D^{-1/2}v$ is equal to the all-ones vector $\one$. Using this we have that
\begin{equation}\label{eq:B2A}
\|I+A\| \ge v^T(I+A)v = v^T D^{-1/2}BD^{-1/2}v = \one^T B \one  =\text{sum}(B).
\end{equation}
This shows the optimal value of~(\ref{eq:Adual}) is at least as great as that of~(\ref{eq:dual}).

Conversely, suppose that $A$ is a feasible solution to~(\ref{eq:Adual}). Since $I+A$ is a nonnegative matrix, by the Perron-Frobenius Theorem it has an entrywise nonnegative eigenvector $u$ (of unit norm) for its maximum eigenvalue. Let $D_u$ be the diagonal matrix with the entries of $u$ on its diagonal, and let $B = D_u(I+A)D_u = (I+A) \circ uu^T$. Since $I+A \succeq 0$, we have that $B \succeq 0$. It is routine to check that $B$ satisfies the other feasibility conditions of~(\ref{eq:dual}). Lastly, note that
\begin{equation}\label{eq:A2B}
\text{sum}(B) = \text{sum}((I+A)\circ uu^T) = u^T (I+A)u = \|I+A\|.
\end{equation}
This shows that $B$ is a feasible dual solution for $\chiv(G)$ with objective value equal to $\|I+A\|$. Therefore, the optimal value of~(\ref{eq:dual}) is at least as great as that of~(\ref{eq:Adual}). Combining this with the above shows that these two optimization programs have the same optimal value.\qeds

The above shows that one can convert a feasible solution to~(\ref{eq:dual}) to a feasible solution to~(\ref{eq:Adual}) with the same or greater objective value, and vice versa. This implies that the conversions take optimal solutions to optimal solutions. Furthermore, we have the following:

\begin{lemma}\label{lem:dualconvert2}
Let $G$ be a nonempty graph. Then there exists an optimal dual solution $B$ for $\chiv(G)$ with strictly positive diagonal if and only if there exists an optimal solution $A$ to~(\ref{eq:Adual}) for $G$ which has a strictly positive maximum eigenvector. Moreover, it can be assumed that $G(B) = G(A)$ and that the multiplicity of $-1$ as an eigenvalue of $A$ is equal to the corank of $B$.
\end{lemma}
\proof
Let  $B$ be an optimal dual solution for $\chiv(G)$ with strictly positive diagonal and set $A =  D^{-1/2}BD^{-1/2} - I$ as in the proof of Lemma~\ref{lem:dualconvert}. Clearly  we have that $G(B) = G(I+A) = G(A)$. Furthermore, as $A$ is optimal to~\eqref{eq:Adual}, Equation~\eqref{eq:B2A} holds throughout with equality, and thus the strictly positive vector $v$ with entries  $v_i = (B_{ii})^{1/2}$   is a maximum eigenvector of $A$. Also note that  $B$ and $I+A=D^{-1/2}BD^{-1/2} $  have the same corank. This implies that the multiplicity of $-1$ as an eigenvalue of $A$ (which is its minimum eigenvalue) is  equal to the corank of $B$.

Conversely, let  $A$ be an optimal solution to~(\ref{eq:Adual}) with a strictly positive maximum eigenvector $u$. Let   $B = D_u(I+A)D_u$ as defined in the proof of Lemma~\ref{lem:dualconvert}. By \eqref{eq:A2B}, $B$ is an optimal dual solution for $\chiv(G)$. Furthermore, we have that $G(B) = G(A)$ and that ${\rm corank}(B) ={\rm corank}(I+A)$. Moreover, $B$ has  strictly positive diagonal.\qeds

\begin{remark} 
Note that if $B$ is connected (i.e., $G(B)$ is a connected graph) then it must have strictly positive diagonal since it is positive semidefinite (thus a zero diagonal entry implies a zero row/column). Conversely, if $A$ is connected, then it will have a strictly positive maximum eigenvector by the Perron-Frobenius Theorem.
\end{remark}

The lemmas above allow us (in some cases) to start with strict complementarity for the factors of a categorical product, move to optimal solutions for~(\ref{eq:Adual}) for the factors, use these to build an optimal solution to~(\ref{eq:Adual}) for the product, and then move back to strict complementarity for the product.

We end this section with the following lemma which we will need for the proof of Corollary~\ref{cor:B}. It says that we can always find \emph{almost} optimal solutions to~(\ref{eq:Adual}) that satisfy certain special properties.

\begin{lemma}\label{lem:irredsol}
Let $H$ be a connected nonempty graph and let $\chiv(H) > k$. Then there exists a feasible solution, $A$, to~(\ref{eq:Adual}) for $H$ such that
\begin{enumerate}
  \packed 
  \item $\|I + A\| > k$;
  \item the maximum eigenvalue of $A$ has multiplicity 1;
  \item there is a maximum eigenvector of $A$ with only positive entries;
  \item the minimum eigenvalue of $A$ is $-1$;
  \item the graph $G(A)$ is connected.
\end{enumerate}
\end{lemma}
\proof
Let $A'$ be an optimal solution to~(\ref{eq:Adual}) for $H$. Then we have that $\|I + A'\| > k$. Also note that since $I + A' \succeq 0$, we have that
$\lambda_{\min}(A') \ge -1$. Furthermore, since multiplying $A'$ by a constant greater than $1$ does not cause a violation of the other constraints on $A'$ and only increases the value of $\|I + A'\|$, we have that $\lambda_{\min}(A') = -1$.

Let $A_H$ be the adjacency matrix of $H$. Define
\[A = \alpha(A' + \varepsilon A_H),\]
where $\varepsilon > 0$ and $\alpha$ is chosen to be positive and such that $\lambda_{\min}(A) = -1$. Note that as $\varepsilon$ approaches 0, the parameter $\alpha$ will approach 1.
Since maximum eigenvalue is a continuous function, for sufficiently small $\varepsilon$, we will have $\|I+A\|>k$. Therefore Conditions (1) and (4) are met by this $A$.

Since $A = \alpha(A' + \varepsilon A_H)$, the graph $G(A)$ is simply $H$, which is connected by assumption. So Condition (5) is satisfied. By the Perron-Frobenius Theorem, the maximum eigenvalue of $A$ has multiplicity 1 and this eigenvalue has an eigenvector whose entries are all positive. Therefore, $A$ meets Conditions (2) and (3).\qeds

\subsection{Vector Hedetniemi}\label{subsec:hedetniemi}

We can now use the formulation of $\chiv$ given in~(\ref{eq:Adual}) to prove the vector coloring analog of Hedetniemi's Conjecture.

\begin{theorem}\label{thm:hedetniemi}
For any graphs $G$ and $H$ we have 
\[\chiv(G \times H)    = \min \{ \chiv(G), \chiv(H) \}.\]
\end{theorem}
\proof
Since vector colorings of the factors induce vector colorings of the product, we have that 
\[\chiv(G \times H)   \le  \min \{ \chiv(G), \chiv(H) \}.\]

To see the other inequality, suppose that $\chiv(G) = s$, $\chiv(H) = t$, and $s \le t$. Let $A_G$ and $A_H$ be optimal solutions to~(\ref{eq:Adual}) for $G$ and $H$ respectively. Define $A =  \frac{1}{t - 1} A_G \otimes A_H$. Since $A_G$ and $A_H$ are optimal, their minimum eigenvalue must be $-1$ (see proof of Lemma \ref{lem:irredsol}), and their maximum eigenvalues must be $s-1$ and $t-1$ respectively. Consequently: 
\begin{itemize}
  \packed 
  \item the minimum eigenvalue of $A$ is
  \[\min \left\{ (-1) \cdot \frac {t-1}{t-1}, \frac {s-1}{t-1} \cdot (-1) \right\} = -1;\]
  \item the maximum eigenvalue of $A$ is
  \[\frac {(t-1)(s-1)}{t-1} = s-1.\] 
\end{itemize}
It follows that $I +A \succeq 0$ and $\| I + A \| = s$. It is easy to verify that $A$ satisfies all the other requirements of~\eqref{eq:Adual}, thus $\chiv(G \times H) \ge s$.\qeds

\begin{remark}\label{eq:theoremmod}
Note that in the case where $\chiv(G) < \chiv(H)$, the matrix $A_H$ in the above proof does not necessarily need to be optimal. Instead, it suffices for it to satisfy the following three properties:  feasibility for~\eqref{eq:Adual} for $H$, $\lambda_{\min}(A_H)=-1$, and $\lambda_{\max}(A_H) \ge \lambda_{\max}(A_G)$ (i.e., its objective value is at least as great as that of $A_G$). This fact is used in Corollary~\ref{cor:B}.
   \end{remark}

A similar proof can be used for the strict vector chromatic number. This has been proven in~\cite{sabvshed}; however, the proof presented here is more direct.

\begin{theorem}\label{thm:hedetniemis}
For any graphs $G$ and $H$ we have 
\[\chisv(G \times H)    = \min \{ \chisv(G), \chisv(H) \}.\]
\end{theorem}
\proof
We use the formulation for $\chisv$ that is analogous to~(\ref{eq:Adual}). This appears in~\cite{Lovasz} and is exactly the same as~(\ref{eq:Adual}) except without the nonnegativity constraint on the entries of $A$. The proof is exactly the same as that of Theorem~\ref{thm:hedetniemi}.\qeds

In the next section we will investigate some properties of vector colorings that are crucial to our proof of Result~\ref{res:A} and, consequently, Result~\ref{res:Acor}. We remark that we have already introduced all that we need for the proofs of Results~\ref{res:B} and~\ref{res:Bcor}. However, when we prove these results we will additionally present some necessary conditions on $G$ for it to induce all of the optimal vector colorings of $G \times H$, and these require some of the notions discussed in Section~\ref{sec:techlems}.

\section{Properties of Vector Colorings}\label{sec:techlems}

In Section~\ref{subsec:Atheorem} we prove Result~\ref{res:A} which shows when the vector colorings of two graphs can be used to construct all of the vector colorings of their categorical product. In order to do this, we will need to develop a bit of theory about optimal vector colorings of graphs. We do this here so that the results in Section~\ref{sec:products} can be presented in a succinct manner.

\subsection{Skeletons and neighborliness}\label{subsec:specialprops}

One issue of importance to us will be when the vector coloring inequality, $p_i^Tp_j \le -1$ for $i \sim j$, is satisfied with equality. Given an optimal vector coloring $\p$ of some graph, we will write $i \sim_\p j$ whenever $i \sim j$ and $p_i^Tp_j = -1$. We will say that such edges are \emph{tight} in $\p$ and all other edges are \emph{slack} in $\p$. This allows us to define the following which will be a key notion in this and later sections:

\begin{definition}
Let $G$ be a graph with optimal vector coloring $\p = (p_1, \ldots, p_n)$. We define the graph $G^\p$ to be the spanning subgraph of $G$ that contains all of the edges of $G$ that are tight in $\p$. We further define the \emph{skeleton} of $G$, denoted $G^{\sk}$, to be the spanning subgraph of $G$ containing only the edges that are tight in every optimal vector coloring of $G$. We write $i \sim_{\sk} j$ if $i$ and $j$ are adjacent in $G^{\sk}$. We will use $N^\p(i)$ and $N^{\p}[i]$ to denote the open and closed neighborhoods of $i$ in $G^\p$, and will similarly use $N^{\sk}(i)$ and $N^{\sk}[i]$ for the same in $G^{\sk}$.
\end{definition}

\begin{remark}\label{rem:tightsk}
Let $\p$ be a max-rank vector coloring. For any pair of vertices satisfying $i \sim_\p j$, Condition $(ii)$ from Theorem~\ref{thm:main} becomes $p_i^TRp_j \le 0$. Moreover, for any pair of vertices satisfying $i \sim_{\sk} j$, Condition $(ii)$ can be replaced by $p_i^TRp_j = 0$.
\end{remark}

Obviously, for any optimal vector coloring $\p$ of $G$, the graph $G^\p$ is not empty (unless $G$ is empty), since otherwise $\p$ could not be optimal. It may be less obvious that $G^{\sk}$ is non-empty whenever $G$ is, but in fact we have the following:

\begin{lemma}\label{lem:sksandwich}
Let $G$ be a graph. Then for any optimal vector coloring $\p$ and optimal dual solution $B$, we have
\[G(B) \subseteq G^{\sk} \subseteq G^\p.\]
Moreover, there exists a max-rank vector coloring $\q$ such that $G^{\sk} = G^\q$.
\end{lemma}
\proof
By complementary slackness, we have that $B_{ij} \ne 0 \Rightarrow M_{ij} = -1$ for $M$ being the Gram matrix of any optimal vector coloring. This proves the first containment $G(B) \subseteq G^{\sk}$, and the second containment is obvious.

For the final claim, if $G^{\sk} = G$, then any max-rank vector coloring will do. Otherwise, suppose that $e \in E(G)$ and $e \not\in E(G^{\sk})$. Then by definition of the skeleton of $G$, there exists some optimal vector coloring $\p^e$ of $G$ in which $e$ is not tight. Let $M^e$ be the Gram matrix of this vector coloring and note that this implies that $M^e_{ij} < -1$ for $ij = e$. Define $M^e$ similarly for all $e \in E(G) \setminus E(G^{\sk})$. Now let $\p'$ be the optimal vector coloring whose Gram matrix is given by
\[M = \frac{1}{|E(G) \setminus E(G^{\sk})|} \sum_{e \in E(G) \setminus E(G^{\sk})} M^e.\]
Then we have that $G^{\p'} = G^{\sk}$. Now let $N$ be the Gram matrix of any max-rank vector coloring of $G$. Obviously, $\frac{1}{2}(M+N)$ is the Gram matrix of some max-rank vector coloring $\q$ and $G^\q = G^{\sk}$.\qeds

We now define the notion of neighborliness:

\begin{definition}
Given an optimal vector coloring $\p$ of a graph $G$, we say that a vertex $i \in V(G)$ is \emph{neighborly} in $\p$ if
\[-p_i \in \cone\left(\left\{p_j : j \sim_\p i\right\}\right).\]
We will simply say that $i$ is \emph{neighborly} if it is neighborly in every optimal vector coloring.
We also define $i \to_\p j$ if
\[-p_i = \sum_{\ell : \ell \sim_\p i} \alpha_\ell p_\ell \text{ for some } \alpha_\ell \ge 0 \text{ where } \alpha_j > 0.\]
We write $i \to j$ if $i \to_\p j$ for all optimal vector colorings $\p$ of $G$.

We also define $D^\p(i) = \{j \in V(G): i \to_\p j\}$ and $D(i) = \{j \in V(G): i \to j\}$. Furthermore, we let $D^\p[i] = \{i\} \cup D^\p(i)$ and $D[i] = \{i\} \cup D(i)$.
\end{definition}

\begin{remark}\label{rem:arrow2sk}
Obviously, $i \to j$ (resp.~$i \to_\p j$) is only possible if $i$ is neighborly (resp.~neighborly in $\p$), and it implies that $i \sim_{\sk} j$ (resp.~$i \sim_\p j$) by definition. Note that it is not clear, and in fact we do not know, whether $i \to j$ or $i \to_\p j$ are symmetric relations.
\end{remark}

It will be useful to express neighborliness in terms of convex hulls instead of conical hulls, which we do in the lemma below.

\begin{lemma}\label{lem:cone2conv}
Let $G$ be a graph with optimal vector coloring $\p$. Then a vertex $i \in V(G)$ is neighborly in $\p$ if and only if
\[\zero \in \conv\left(\{p_j : j \in N^\p[i]\} \right).\]
Moreover, if $\zero = \sum_{j \in N^\p[i]} \alpha_j p_j$, then $\sum_{j \in N^\p[i]} \alpha_j  = \chiv(G)\alpha_i$.
\end{lemma}
\proof
Let $t = \chiv(G)$ and suppose that $\zero = \alpha_i p_i + \sum_{j \sim_\p i} \alpha_j p_j$. Taking inner product with $p_i$ on both sides reveals $0 = \alpha_i (t-1) - \sum_{j \sim_\p i} \alpha_j$, and thus $\alpha_i + \sum_{j \sim_p i} \alpha_j = \alpha_i t$. Thus we have proven the second claim.

The above shows that if $\zero = \sum_{j \in N^\p[i]} \alpha_j p_j$ where the righthand side is a convex combination, then the coefficient $\alpha_i$ of $p_i$ is nonzero. Thus, if $\zero \in \conv\left(\{p_j : j \in N^\p[i]\} \right)$, then $-p_i \in \cone\left(\{p_j : j \in N^\p(i)\} \right)$, i.e., $i$ is neighborly in $\p$. The other direction holds since $-p_i \in \cone\left(\{p_j : j \sim_\p i\} \right)$ implies that we can find nonnegative coefficients $\alpha_j$ such that $-p_i = \sum_{j \sim_p i} \alpha_j p_j$ and thus $\zero = p_i + \sum_{j \sim_p i} \alpha_j p_j$. Rescaling the righthand side gives a convex combination equal to $\zero$.\qeds

Both the convex hull and conical hull perspectives are useful. The convex hull view is used to prove Lemma~\ref{lem:nbrlysk}, whereas conical hulls are essential in Lemma~\ref{lem:conespan}.

It turns out that for a vertex to be neighborly, it suffices for it to be neighborly in some max-rank vector coloring, and moreover the conical/convex combination witnessing neighborliness can be fixed for all vector colorings:

\begin{lemma}\label{lem:nbrlymr}
Let $G$ be a graph with max-rank vector coloring $\p$. Then, $i\to_\p j$ implies that $i\to j$.   Furthermore, $i \in V(G)$ is neighborly if and only if it is neighborly in $\p$.
\end{lemma}
\proof
Let $\p$ be a max-rank vector coloring and suppose that $i$ is neighborly in $\p$, i.e., that $\zero = \sum_{j \in N^\p[i]} \alpha_j p_j$ for some $\alpha_j \ge 0$ for all $j$. Now let $\alpha$ be the vector of coefficients from the righthand side, but extended to $|V(G)|$ coordinates by adding zeros in the appropriate places. Also, let $P$ be the matrix whose rows are the $p_\ell^T$ for $\ell \in V(G)$. Then we can rewrite the equation above as $\zero = P^T\alpha$. Of course, this implies that $PP^T\alpha = \zero$ and thus $\alpha$ is a vector in the kernel of the Gram matrix of the vector coloring $\p$. Now let $\q$ be some other optimal vector coloring of $G$ and let $Q$ be the matrix whose rows are the $q_\ell^T$. Then by Lemma~\ref{lem:image} we have that $\ker(PP^T) \subseteq \ker(QQ^T)$. Therefore we have that $QQ^T\alpha = \zero$ and this is equivalent to $Q^T\alpha = \zero$. Of course, the latter is equivalent to $\sum_{j \in N^\p[i]} \alpha_j q_j = \zero$. Note that we are not done yet because we are still summing over $j \in N^\p[i]$. However, suppose that $j' \not\sim_\q i$ for some $j'$ such that $\alpha_{j'} > 0$. Then $q_i^Tq_{j'} < -1$, and thus taking inner product with $q_i$ on both sides of the above equation gives
\begin{align*}
0 &= \alpha_iq_i^T q_i + \sum_{j \sim_p i} \alpha_j q_i^Tq_j \\
&< \alpha(t-1) - \sum_{j \sim_\p i} \alpha_j \\
&= p_i^T(\alpha_i p_i + \sum_{j \sim_\p i} \alpha_j p_j) = 0,
\end{align*}
a clear contradiction. Thus we can conclude that $j \sim_\q i$ for all $j$ such that $\alpha_j > 0$, and so $\sum_{j \in N^\q[i]} \alpha_j q_j = \zero$. This shows that $i$ is neighborly in $\q$ and that $i \to_\q j$ for all $j$ such that $\alpha_j > 0$. Since $\q$ was an arbitrary optimal vector coloring, we have that $i$ is neighborly and that $i \to j$ for all $j$ such that $\alpha_j > 0$. If $i \to_\p j$, then by definition we could have chosen our convex combination such that $\alpha_j > 0$. Therefore, if $i \to_\p j$, then $i \to_\q j$ for all optimal vector colorings $\q$ and thus $i \to j$.\qeds

\begin{remark}\label{rem:dualnbrly}
The coefficients in these convex combinations are playing the role of the rows/columns of an optimal dual solution. Indeed, if $B$ is any optimal dual solution and $P$ is a matrix whose rows are the vectors of an optimal vector coloring of $G$, then we have $PP^TB = 0$ by complementary slackness, and therefore $P^TB = 0$. This latter equation is equivalent to $\sum_{j}B_{ji} p_j = \zero$ for all $i \in V(G)$. Some of the rows/columns of $B$ may be zero, but if the $i^\text{th}$ column is nonzero, then this equation shows that $i$ is neighborly (since it must hold for all optimal vector colorings by complementary slackness). Note that since $B$ is positive semidefinite, its $i^\text{th}$ row/column being nonzero is equivalent to its $i^\text{th}$ diagonal entry being nonzero. Thus an optimal dual solution $B$ with $B_{ii} > 0$ implies that vertex $i$ is neighborly.

The above is one of the reasons why simply proving Results~\ref{res:Bcor} and~\ref{res:Acor} directly would be easier: we could use properties of the type of dual solution which is assumed to exist in those results in order to obtain properties of the vector colorings of the graph(s) in question. This is quicker than building up theory about optimal vector colorings as we are doing here, but it would not allow us to prove the necessary and sufficient conditions of Results~\ref{res:B} and~\ref{res:A}. We note that we do not know how to go in the other direction: to use the convex combinations witnessing neighborliness to construct an optimal dual solution.
\end{remark}

\begin{lemma}\label{lem:conespan}
Let $G$ be a graph with optimal vector coloring $\p$. If every vertex of $G$ is neighborly in $\p$, then
\[\cone\left(\left\{p_i - p_j : i,j \in V(G), \ i \to_\p j\right\}\right) = \cone \left(\left\{p_i : i \in V(G)\right\}\right) = \spn(\p).\]
If every vertex of $G$ is neighborly, then also
\[\cone\left(\left\{p_i - p_j : i,j \in V(G), \ i \to j\right\}\right) = \spn(\p).\]
\end{lemma}
\proof
We will show that
\[\cone\left(\left\{p_i - p_j : i,j \in V(G), \ i \to_\p j\right\}\right) \supseteq \cone \left(\left\{p_i : i \in V(G)\right\}\right) \supseteq \spn(\p)\]
which proves the first claim since it is obvious that both cones are contained in $\spn(\p)$.

To show the first containment we only need to show that $p_i \in \cone(\{p_i - p_j: i \to_\p j\})$ for all $i \in V(G)$. By the assumption of neighborliness we have that $-p_i = \sum_{j \in D^\p(i)} \alpha_j p_j$ where $\alpha_j \ge 0$ for all $j$. We can rewrite this as $p_i = \sum_{j \in D^\p(i)} \alpha_j (-p_j)$. This implies that
\[\left(1 + \sum_{j \in D^\p(i)} \alpha_j\right)p_i = \sum_{j \in D^\p(i)} \alpha_j (p_i - p_j).\]
Since the coefficient on the lefthand side is strictly positive, this shows that $p_i \in \cone(\{p_i - p_j: i \to_\p j\})$ as desired. Thus we have proven the first containment.

To show the second containment, note that by assumption of neighborliness we have that $-p_i \in \cone(\{p_j : j \sim_\p i\}) \subseteq \cone(\{p_j : j \in V(G)\})$ for all $i \in V(G)$. This already implies that $\cone(\{p_j : j \in V(G)\}) = \spn(\p)$ and so we are done with the first claim. The second claim follows exactly as the first.\qeds

The next two results concern the skeleton of a graph. The first one relates this notion to that of neighborliness.

\begin{lemma}\label{lem:nbrlysk}
Let $G$ be a nonempty graph. Then $i \in V(G)$ is neighborly if and only if it is not isolated in $G^{\sk}$.
\end{lemma}
\proof
Let $\q$ be an optimal vector coloring of $G$ such that $G^\q = G^{\sk}$. If $i \in V(G)$ is isolated in $G^{\sk}$, then it is isolated in $G^\q$ and so obviously $-q_i \not\in \cone(\{q_j : j \in N^\q(i)\}) = \cone(\varnothing)$ and thus $i$ is not neighborly.

Conversely, suppose that $i$ is not neighborly. Then by definition $i$ is not neighborly in some optimal vector coloring $\p$ of $G$. Therefore, by Lemma~\ref{lem:cone2conv}, we have that $\zero \not\in \conv(\{p_j : j \in N^\p[i]\})$. Since $\conv(\{p_j : j \in N^\p[i]\})$ is a compact convex set, by the Hyperplane Separation Theorem there exists a vector $v$ such that $v^Tw < 0$ for all $w \in \conv(\{p_j : j \in N^\p[i]\})$ and $v^T\zero = 0$. We can further choose $v$ (by rescaling if necessary) so that $v^Tp_j < -1$ for all $j \in N^\p[i]$. We will use $v$ to show that we can replace $p_i$ with some $p'_i$ such that $\|p'_i\|^2 = t-1$ and $p_j^Tp'_i < -1$  for all $j \sim i$. This will show that $i$ must be isolated in $G^{\sk}$.

Consider the convex combination $p_\varepsilon = (1 - \varepsilon) p_i + \varepsilon v$ for $0 < \varepsilon < 1$. For $j \sim_\p i$, it is easy to see that
\[p_\varepsilon^Tp_j = (1 - \varepsilon) p_i^Tp_j + \varepsilon v^Tp_j < (1 - \varepsilon) (-1) + \varepsilon (-1) = -1.\]
Since $p_i^Tp_j < -1$ for all $j \sim i$ such that $j \not\sim_\p i$, we can pick $\varepsilon$ close enough to 0 so that
\[p_\varepsilon^Tp_j = (1-\varepsilon) p_i^Tp_j + \varepsilon v^Tp_j < -1\]
for all such $j$. Thus for $\varepsilon$ sufficiently small, $p_\varepsilon$ satisfies $p_\varepsilon^Tp_j < -1$ for all $j \sim i$. To finish, we must show that for sufficiently small $\varepsilon$, the vector $p_\varepsilon$ has norm at most that of $p_i$. If this is true, then we can rescale $p_\varepsilon$ so that it has norm squared equal to $t - 1$ while still maintaining these strict inequalities. Since $v^Tp_i < 0$, the vectors $v$ and $p_i$ form an obtuse angle at the origin, and so it is ``geometrically obvious" that for small enough $\varepsilon$ the vector $p_\varepsilon$ has strictly smaller norm than $p_i$. However, we will give a rigorous proof.

Let $s = \|v\|^2$. We have that
\[\|p_\varepsilon\|^2 = (1 - \varepsilon)^2 (t-1) + 2\varepsilon(1 - \varepsilon)v^Tp_i + \varepsilon^2 s < (t-1) + \varepsilon (2v^Tp_i + \varepsilon(s - 2v^Tp_i)).\]
Since $v^Tp_i < 0$, for sufficiently small $\varepsilon$ the $2v^Tp_i + \varepsilon(s - 2v^Tp_i)$ term is strictly negative, and so we are done.

So if we replace $p_i$ with a rescaled version of $p_\varepsilon$ for sufficiently small $\varepsilon$, we will obtain an optimal vector coloring of $G$ in which every edge incident to $i$ is slack. This implies that $i$ must be isolated in $G^{\sk}$ as desired.\qeds

The next lemma relates properties of a graph to those of its skeleton, specifically their vector chromatic number and rank.

\begin{lemma}\label{lem:skcomps}
Let $G$ be a graph and let $G_\ell$ for $\ell = 1, \ldots, m$ be the connected components of $G^{\sk}$ that are not isolated vertices, and let $S$ be the set of isolated vertices of $G^{\sk}$. Then $\chiv(G_\ell) = \chiv(G)$ for all $\ell \in [m]$. Furthermore,
$\rk(G) = |S| + \sum_{\ell=1}^m \rk(G_\ell) = \rk(G^{\sk})$ and $(G^{\sk})^{\sk} = G^{\sk}$.
\end{lemma}
\proof
Let $t = \chiv(G)$. First note that $\chiv(G_\ell) \le t$ since $G_\ell$ is a subgraph of $G$. Now let $\p$ be an optimal vector coloring of $G$ such that $G^\p = G^{\sk}$. If $M^\p$ is the Gram matrix of $\p$, then we have that $M^\p_{ij} < -1$ for all $i,j \in V(G)$ such that $i \sim j$ and $i \not\sim_{\sk} j$. Now suppose that $k \in [m]$ is such that $\chiv(G_{k}) < t$. For each $\ell \ne k$, let $\q^\ell$ be a vector $t$-coloring of $G_\ell$. Let $\q^{k}$ be an optimal vector coloring of $G_{k}$ that has been globally rescaled so that $\|q^{k}_i\|^2 = t-1$ for all $i \in V(G_k)$. Note that this implies that $(q^{k}_i)^Tq^{k}_j < -1$ for all $i \sim j$ in $G_{k}$. Let $M^\ell$ be the Gram matrix of the vectors in $\q^\ell$ for each $\ell \in [m]$. Define $M$ to be the block diagonal matrix with blocks given by the $M^\ell$ for $\ell \in [m]$, and additionally $M_{ii} = t-1$ for all $i \in S$. Note that $M_{ii} = t-1$ for all $i \in V(G)$ and $M_{ij} \le -1$ for $i \sim j$ unless $i \not\sim_{\sk} j$, in which case $M^\p_{ij} < -1$. Thus it is easy to see that for sufficiently small $\varepsilon > 0$, the convex combination $(1 - \varepsilon)M^\p + \varepsilon M$ is the Gram matrix of an optimal vector coloring of $G$ such that every edge in $G_{k}$ is slack. This is a contradiction since all of the edges in $G_{k}$ are contained in $G^{\sk}$. This proves the first claim.

To show that $\rk(G) = |S| + \sum_{\ell=1}^m \rk(G_\ell)$, note that the first claim implies that $\dim \spn(\{q_i : i \in V(G_\ell)\}) \le \rk(G_\ell)$ for all $\ell \in [m]$. Moreover, $\dim \spn(\{q_i : i \in S\}) \le |S|$ trivially holds. Thus we have that for any optimal vector coloring $\q$ of $G$,
\[\dim \spn(\q) \le |S| + \sum_{\ell = 1}^m \dim \spn(\{q_i : i \in V(G_\ell)\}) \le |S| + \sum_{\ell = 1}^m \rk(G_\ell).\]
This proves that $\rk(G) \le |S| + \sum_{\ell=1}^m \rk(G_\ell)$. The proof of the other inequality is similar to the proof of the first claim above. For each $\ell \in [m]$, we let $M^\ell$ be the Gram matrix of a \emph{max-rank} vector coloring of $G_\ell$. Thus $\rk(M^\ell) = \rk(G_\ell)$ for all $\ell \in [m]$. Now let $M$ be defined as above, as the block diagonal matrix with blocks given by the $M^\ell$'s, and with $1 \times 1$ blocks consisting of a single $t-1$ entry for each $i \in S$. Then we have $\rk(M) = |S| + \sum_\ell \rk(M^\ell) = |S| + \sum_\ell \rk(G_\ell)$. As above, for sufficiently small $\varepsilon > 0$, the convex combination $(1 - \varepsilon)M^\p + \varepsilon M$ is the Gram matrix of an optimal vector coloring of $G$. Furthermore, this convex combination has rank at least that of $M$, and so we have shown $\rk(G) = \sum_{\ell=1}^m \rk(G_\ell)$. The next equality in the lemma follows from $(G^{\sk})^{\sk} = G^{\sk}$ whose proof should be clear at this point.\qeds

Note that the above lemma implies that $\rk(G \cup H) = \rk(G) + \rk(H)$ whenever $\chiv(G) = \chiv(H)$.

\subsection{Examples of skeletons}\label{subsec:skeletons}

In the previous section we investigated several properties of the skeleton of a graph. However, we have not yet seen any actual examples of these objects. Perhaps the skeleton of a graph is always just the graph itself? Here we will determine the skeletons of some basic graphs, as well as some more complicated examples.

It is not hard to see that the skeleton of a complete graph is itself: indeed a complete graph has a unique vector coloring and in this vector coloring all of the edges are tight (see the discussion following Corollary~\ref{cor:B2} for a quick proof of this). This is a special case of edge-transitive graphs, graphs such that for any two edges there is an automorphism mapping the first edge to the second. Since all of the edges of such a graph are ``the same", they are either all in the skeleton or none are. By Lemma~\ref{lem:skcomps}, the latter is impossible (unless the graph is empty). Another example of graphs that are equal to their skeletons are bipartite graphs. These graphs have vector chromatic number equal to 2, and it is not difficult to see that in any vector 2-coloring the vectors assigned to the ends of any edge must have the form $v, -v$ for some \emph{unit} vector $v$, and thus their inner product is $-1$, i.e., the edge is tight.

The smallest example of a graph which is not equal to its skeleton is $K_3$ plus a vertex adjacent to one of the vertices of the $K_3$. It is easy to see that this graph has an optimal vector coloring in which the edge incident to the degree one vertex is not tight, thus the skeleton of this graph is $K_3$ plus an isolated vertex. We can change $K_3$ to $K_n$ for $n \ge 3$, and change the single edge to a longer path to construct similar examples where the skeleton is now $K_n$ plus some number of isolated vertices. We can put another $K_n$ at the other end of the path to obtain an example of a connected graph whose skeleton has  more than one nontrivial connected component. Shortening the path back to a single edge, we obtain connected graph whose skeleton is two $K_n$'s. We can play around with this and similar constructions to obtain examples of graphs and skeletons which allow us to build up some intuition about this notion. In fact, already the first example was significant to our intuition during the development of this work.

A slightly more advanced example comes from the graph $H_{n,k}$ investigated in~\cite{UVC2}. This graph has the even weight binary strings of length $n$ as its vertices, two being adjacent if they are at Hamming distance exactly $k$ (where $k$ is restricted to being even). We showed in~\cite{UVC2} that these graphs are uniquely vector colorable whenever $n \le 2k-2$, and moreover they remain so (with the same unique vector coloring) if any number of edges are added between vertices at Hamming distance greater than $k$. Finally, these added edges will always be slack in the unique vector coloring of this graph. This gives a large family of graphs whose skeleton is $H_{n,k}$.

\section{Vector Colorings of the Categorical Product}\label{sec:products}

In this section we will prove Results~\ref{res:B}--\ref{res:Acor} showing when the optimal vector colorings of a product are determined by the vector colorings of the factors. We will begin with the case where one factor has strictly smaller vector chromatic number, i.e., Results~\ref{res:B} and~\ref{res:Bcor}. These will be relatively straightforward, but Results~\ref{res:A} and~\ref{res:Acor}, presented in Section~\ref{subsec:Atheorem}, will require more substantial proofs. In Section~\ref{subsec:1walkreg} we will consider 1-walk-regular graphs, showing that they always satisfy strict complementarity. Finally, in Section~\ref{subsec:implications} we prove a vector coloring analog of the Duffus, Sands, and Woodrow result that $(A_n) \Rightarrow (B_n)$ for all $n \in \mathbb{N}$.

\subsection{Factors with different vector chromatic numbers}\label{subsec:Btheorem}

In the case where $\chiv(G) < \chiv(H)$, only the vector colorings induced by $G$ are optimal, not those induced by $H$. The maximum rank of an optimal vector coloring induced by $G$ is $\rk(G)$, and so $\rk(G \times H) = \rk(G)$ is obviously necessary for all of the optimal vector colorings of $G \times H$ to be induced by $G$. Below we show that this is also sufficient.

\begin{theorem}\label{thm:B}
Let $G$ and $H$ be graphs such that $\chiv(G) < \chiv(H)$. Then $\rk(G \times H) \ge \rk(G)$ and equality holds if and only if every optimal vector coloring of $G \times H$ is induced by $G$.
\end{theorem}
\proof
By Theorem~\ref{thm:hedetniemi}, the vector colorings of $G \times H$ induced by the optimal vector colorings of $G$ are in fact optimal for $G \times H$. Moreover, any such induced vector coloring of $G \times H$ spans the same dimension as the corresponding vector coloring of $G$, since it uses exactly the same set of vectors. Therefore $\rk(G \times H) \ge \rk(G)$ and $\rk(G \times H) > \rk(G)$ is only possible if there is some optimal vector coloring of $G \times H$ that is not induced by $G$. Thus we have proven one direction of the claim.

Now suppose that $\rk(G \times H) = \rk(G)$, and let $\p$ be an optimal vector coloring of $G$ that spans $\mathbb{R}^d$ for $d = \rk(G)$. Let $\q$ be the vector coloring of $G \times H$ induced by $\p$, i.e., $q_{i\ell} = p_i$ for all $i \in V(G)$, $\ell \in V(H)$. Let $P$ be the matrix whose rows are the $p_i$ and note that the matrix whose rows are the $q_{i\ell}$ is $Q = P \otimes \one$. Also, let $I_1$ denote the $1 \times 1$ identity. Since $\q$ is a max-rank vector coloring of $G \times H$ by assumption, from Theorem~\ref{thm:main} we have that the Gram matrix of any optimal vector coloring of $G \times H$ is equal to
\[Q(I+R)Q^T = \left(P \otimes \one\right)\left((I+R) \otimes I_1\right)\left(P \otimes \one\right)^T = P(I+R)P^T \otimes J,\]
where $R$ is a symmetric matrix satisfying
\begin{align*}
&q_{i\ell}^T R q_{i\ell} = 0 \text{ for all } i,\ell; \\
&q_{i\ell}^T R q_{jk} \le -1 - q_{i\ell}^Tq_{jk} \text{ for } (i,\ell) \sim (j,k); \\
&I+R \succeq 0.
\end{align*}

But now we see that we are already done. Indeed, the $(i,\ell)(j,k)$-entry of the matrix $P(I+R)P^T \otimes J$ is merely the $ij$-entry of $P(I+R)P^T$, which clearly only depends on $i$ and $j$. Therefore, any such vector coloring is induced by $G$.\qeds

As a corollary we have the following:

\begin{corollary}\label{cor:B}
Let $G$ and $H$ be graphs such that $\chiv(G) < \chiv(H)$ and $H$ is connected. If $G$ admits a strictly complementary dual solution with strictly positive diagonal, then every optimal vector coloring of $G \times H$ is induced by $G$.
\end{corollary}
\proof
Note that the claim holds trivially if $G$ is empty. Thus we may assume that $G$ is nonempty, and $H$ must then be nonempty as $\chiv(H) > \chiv(G)$. We now show that under the hypotheses of the corollary we have that $\rk(G \times H) = \rk(G)$ and then we can apply Theorem~\ref{thm:B}. As the inequality $\rk(G \times H) \ge \rk(G)$ is always true, it remains to show the reverse inequality, i.e., $\rk(G \times H) \le  \rk(G)$. For this, it suffices to find an optimal dual solution $B'$ for $\chiv(G \times H)$ which has corank equal to $\rk(G)$. Indeed, in this case we have that
\begin{equation}\label{eq:cdver}
\rk(G\times H)\le {\rm corank}(B')=\rk(G)\le \rk(G\times H),
\end{equation}
and thus we have equality throughout in \eqref{eq:cdver}. To show the existence of a matrix $B'$ with these properties, by Lemma~\ref{lem:dualconvert2} it suffices to find an optimal solution $A$ to~(\ref{eq:Adual}) for $G\times H$ which has $-1$ as an eigenvalue with multiplicity $\rk(G)$ and an entrywise positive maximum eigenvector.

Let $B$ be a strictly complementary dual solution for $\chiv(G)$ with strictly positive diagonal. Then $\text{corank}(B) = \rk(G)$. By Lemma~\ref{lem:dualconvert2}, there exists an optimal solution $A_G$ to~(\ref{eq:Adual}) for $G$ which has a strictly positive maximum eigenvector. Moreover, the multiplicity of $-1$ as an eigenvalue of $A_G$ is equal to $\text{corank}(B) = \rk(G)$.

Let $A_H$ be a solution to~(\ref{eq:Adual}) for $H$ of value strictly greater than $\chiv(G)$ with the additional properties guaranteed by Lemma~\ref{lem:irredsol}. Also, let $\lambda$ and $\mu$ be the maximum eigenvalues of $A_G$ and $A_H$ respectively. By optimality, both $A_G$ and $A_H$ have least eigenvalue $-1$, and as $\chiv(G) < \chiv(H)$ we have that $\lambda < \mu$. By the proof of Theorem~\ref{thm:hedetniemi} (recall Remark~\ref{eq:theoremmod}), the matrix $A = \frac{1}{\mu}A_G \otimes A_H$ is an optimal solution to~(\ref{eq:Adual}) for $G \times H$.

Next, we show that $\text{corank}(I+A) = \rk(G)$. As $G$ is nonempty, we have that $\chiv(H) > \chiv(G) \ge 2$, and so $\mu > \lambda \ge 1$. Thus, the minimum eigenvalue of $A_G \otimes A_H$ is equal to $-\mu$. Furthermore, as $\mu$ is a simple eigenvalue of $A_H$ (Lemma~\ref{lem:irredsol}) and the multiplicity of  $-1$ as an eigenvalue of $A_G$ is $\rk(G)$ (Lemma \ref{lem:dualconvert2}), the multiplicity of $-1$ as an eigenvalue of $A$ is $\rk(G)$.

Lastly, we show that $A$ has a positive maximum eigenvector. Note that the maximum eigenvalue of $A$ is equal to $\lambda$. By Lemma~\ref{lem:irredsol}, $A_H$ has a strictly positive maximum eigenvector. Additionally, by Lemma~\ref{lem:dualconvert2}, $A_G$ also has a strictly positive maximum eigenvector. Taking the Kronecker product of these two eigenvectors we obtain a strictly positive maximum eigenvector of $A$.\qeds

\begin{remark}\label{rem:connected}
We do not really lose anything by assuming that $H$ is connected in the above corollary. Indeed, $H$ must be connected for every optimal vector coloring of $G \times H$ to be induced by $G$. To see this, suppose that $H$ has connected components $H_1, \ldots, H_k$ for $k \ge 2$. Consider any optimal vector coloring of $G \times H$ and note that we can obtain a new optimal vector coloring by applying an arbitrary orthogonal transformation to all of the vectors assigned to vertices in $V(G) \times V(H_1)$, and fixing the other vectors. It is easy to see that the vector assigned to $(i,\ell)$ for $\ell \in V(H_1)$ in this new coloring will be different than that assigned to $(i,\ell')$ for $\ell' \notin V(H_1)$. But this is not possible for a vector coloring induced by $G$.
\end{remark}

We immediately obtain the following corollary when $G$ is uniquely vector colorable:

\begin{corollary}\label{cor:B2}
Let $G$ be a uniquely vector colorable graph for which there exists a strictly complementary dual solution with strictly positive diagonal. If $H$ is connected and $\chiv(G) < \chiv(H)$, then $G \times H$ is uniquely vector colorable.\qeds
\end{corollary}

Corollary~\ref{cor:B2} above generalizes a result of Pak and Vilenchik~\cite{Pak}. They show that if an $r$-regular graph $H$ with eigenvalues $\lambda_1 \ge \lambda_2 \ge \ldots \ge \lambda_n$ satisfies $\lambda(H) < r/(m-1)$, where $\lambda(H) = \max_{i \ge 2} |\lambda_i |$, then the product $K_m \times H$ is uniquely vector $m$-colorable. It was in fact this result that originally inspired our Theorem~\ref{thm:B}. It is not immediately obvious why this result is implied by Corollary~\ref{cor:B2}, so we give a brief explanation.

Suppose $H$ is as described above. First note that $K_m$ has a unique optimal vector coloring with Gram matrix equal to $mI - J$, and has $\frac{1}{m}J$ as a strictly complementary dual solution. To see that this vector coloring is unique, note that the Gram matrix $M$ of any other vector $m$-coloring of $K_m$ would necessarily have $\text{sum}(M) < \text{sum}(mI-J) = 0$, and thus $\one^TM\one < 0$, a contradiction. Next, $\lambda(H) < r/(m-1)$ implies that $\lambda_2 \ne \lambda_1 = r$ and so $r$ is a simple eigenvalue. Since $H$ is regular, this implies that $H$ is connected. Also,
\[-\lambda_n = |\lambda_n| \le \lambda(H) < \frac{r}{m-1} \ \Longrightarrow \ 1 - \frac{r}{\lambda_n} > m.\]
However, $1 - r/\lambda_n$ is actually a lower bound on the vector chromatic number of $H$. In fact, the vector chromatic number of any graph $G$ is equal to the maximum of $1 - \lambda_{\max}(A)/\lambda_{\min}(A)$ where $A$ ranges over nonnegative symmetric matrices such that $A_{ij} = 0$ if $i \not\sim j$. This formulation for $\chi_v$ can be easily derived from~(\ref{eq:Adual}). Therefore, if $H$ satisfies the hypotheses of the Pak and Vilenchik result, then $H$ is connected and $\chiv(H) > m$, which means that it satisfies the hypotheses of Corollary~\ref{cor:B2} above.

Pak and Vilenchik also note that their result implies that $K_m \times H$ is uniquely $m$-colorable as well. However, it is known~\cite{greenwell} that if $H$ is connected and $\chi(H) > m$, then $K_m \times H$ is uniquely $m$-colorable. Since their hypotheses imply that $\chiv(H) > m$, they also imply that $\chi(H) > m$, and so the classical result is already more general in this regard.

Corollary~\ref{cor:B2} allows one to build many examples of uniquely vector colorable graphs. One could take $G$ to be any of the Kneser or $q$-Kneser graphs, which were proven to be uniquely vector colorable in~\cite{UVC1}. One could also let $G$ be one of the Hamming graphs proven to be uniquely vector colorable in~\cite{UVC2}. As long as $H$ is connected and has strictly larger vector chromatic number than $G$, then $G \times H$ is uniquely vector colorable.

Note that the assumption that $B$ had strictly positive diagonal was essential to our proof of Corollary~\ref{cor:B}. This is because if $B$ did not have this property, and we performed the same steps to obtain $B'$, then $\text{corank}(B')$ would be greater than $\rk(G)$. This does not prove that it is necessary for $G$ to have a strictly complementary dual solution with strictly positive diagonal in order for the conclusion of Corollary~\ref{cor:B} to hold; in fact we do not even know if $G$ necessarily must have a strictly complementary dual solution at all. However, we can show that the analogous condition for neighborliness (that every vertex is neighborly, recall Remark~\ref{rem:dualnbrly}) is necessary:

\begin{lemma}\label{lem:Bnbrly}
Suppose that $G$ and $H$ are graphs such that $\chiv(G) < \chiv(H)$ and $H$ is connected. If every optimal vector coloring of $G \times H$ is induced by $G$, then every vertex of $G$ is neighborly, i.e., $G^{\sk}$ has no isolated vertices.
\end{lemma}
\proof
Suppose vertex $i^* \in V(G)$ is not neighborly. Since $i^*$ is not neighborly, it is an isolated vertex in $G^{\sk}$ by Lemma~\ref{lem:nbrlysk}. Thus there exists an optimal vector coloring $\p$ of $G$ such that $i^*$ is isolated in $G^\p$, i.e., such that $p_{i^*}^Tp_j < -1$ for all $j \sim i^*$. Let $\q$ be the vector coloring of $G \times H$ induced by $\p$, so $q_{i\ell} = p_i$ for all $i \in V(G)$, $\ell \in V(H)$. Now fix some $\ell^* \in V(H)$ and note that $q_{i^* \ell^*}^Tq_{jk} < -1$ for all $j \sim i^*$, $k \sim \ell^*$. It is easy to see that applying some small rotation to $q_{i^* \ell^*}$ and fixing all other vectors in $\q$ will not break any of the properties required of an optimal vector coloring of $G \times H$, but the new vector coloring will not be induced by $G$, since the vector assigned to $(i^*, \ell^*)$ will not be the same as that assigned to $(i^*, k)$ for $k \ne \ell^*$. This is a contradiction to our assumption that every optimal vector coloring of $G \times H$ is induced by $G$, and so every vertex of $G$ must be neighborly.\qeds

Putting all of the above together we obtain the following:

\begin{theorem}\label{thm:Bsandwich}
Let $G$ be a nonempty graph. Then the following statements satisfy $(1) \Rightarrow (2) \Rightarrow (3) \Rightarrow (4)$:
\begin{enumerate}
\item $G$ admits a strictly complementary dual solution with strictly positive diagonal.
\item For any connected graph $H$ with $\chiv(H) > \chiv(G)$, every optimal vector coloring of $G \times H$ is induced by $G$.
\item There exists a connected graph $H$ with $\chiv(H) > \chiv(G)$ such that every optimal vector coloring of $G \times H$ is induced by $G$.
\item Every vertex of $G$ is neighborly.
\end{enumerate}
\end{theorem}

It is an interesting question whether any of these implications can be reversed. Or, even better, if one can show that $(4)$ implies $(1)$. This is related to whether the dependencies arising from neighborly vertices can be used to construct optimal dual solutions for vector chromatic number.

\subsection{Factors with the same vector chromatic number}\label{subsec:Atheorem}

In the case where both factors have the same vector chromatic number, each factor can induce optimal vector colorings of the product. The Gram matrix of such an induced vector coloring has the form $M \otimes J$ or $J \otimes N$ where $M$ and $N$ are the Gram matrices of some optimal vector colorings of the first and second factor respectively. Such vector colorings are easy to recognize: they are induced by $G$ if the $(i,\ell)(j,k)$-entry of the Gram matrix depends only on $i$ and $j$, and they are induced by $H$ if it only depends on $\ell$ and $k$. However, one can also take convex combinations of vector colorings induced by each of the factors. This results in a vector coloring whose Gram matrix has the form $\alpha M \otimes J + \beta J \otimes N$ where $0 \le \alpha = 1 - \beta \le 1$. Such a mixing of vector colorings induced by the factors is more difficult to recognize, and indeed we do not know a simple necessary and sufficient condition for when a vector coloring has this form. However, we are able to prove the following:

\begin{lemma}\label{lem:convcombvc}
Let $G$ and $H$ be graphs with $\chiv(G) = \chiv(H)$. Suppose that $M \otimes J + J \otimes N$ is an optimal primal solution for $\chiv(G \times H)$ where $M, N \succeq 0$ are matrices indexed by $V(G)$ and $V(H)$ respectively. Then $M \otimes J + J \otimes N$ is a convex combination of optimal vector colorings of $G \times H$ induced by $G$ and $H$.
\end{lemma}
\proof
If either factor is empty then both must be and the product will be as well. In this case we must have that $M$ and $N$ are both zero matrices and we are done. So we may suppose that $G$ and $H$ are not empty.

Let $t = \chiv(G) = \chiv(H) = \chiv(G \times H)$. Considering the diagonal, we see that for any $i \in V(G)$ and $\ell \in V(H)$,
\[t - 1 = (M \otimes J + J \otimes N)_{i\ell,i \ell} = M_{ii} + N_{\ell \ell}.\]
Fixing $i$ and varying $\ell$, and vice versa, shows that there exists some $\gamma \in \mathbb{R}$ such that
\[M_{ii} = \gamma \quad \text{and} \quad N_{\ell\ell} = t-1 - \gamma\]
for all $i \in V(G), \ell \in V(H)$. Since $M, N \succeq 0$, we see that $\gamma \ge 0$ and $t - 1 - \gamma \ge 0$. Furthermore, if either constant is equal to zero, then one of $M$ and $N$ must be zero and then we are in the case where the vector coloring is induced by a single factor. So we may assume that $0 < \gamma < t-1$. Let $\alpha = \gamma/(t-1)$ and $\beta = (t-1-\gamma)/(t-1)$, and note that $\alpha,\beta > 0$ and $\alpha + \beta = 1$. We will show that $\alpha^{-1} M$ and $\beta^{-1}N$ are Gram matrices of optimal vector colorings of $G$ and $H$ respectively. Note that they are both positive semidefinite and have constant diagonal equal to $t-1$ by construction. So we only need to check that their entries corresponding to edges are at most $-1$.

Suppose for contradiction that $\alpha^{-1} M_{ij} > -1$, for some $i \sim j$. Then for any $\ell \sim k$ in $H$, we have that $(i,\ell) \sim (j,k)$ in $G \times H$, and so
\[-1 \ge (M \otimes J + J \otimes N)_{i\ell, jk} = \alpha(\alpha^{-1} M_{ij}) + \beta(\beta^{-1}N_{\ell k})\]
As the righthand side is a convex combination of $\alpha^{-1}M_{ij}$ and $\beta^{-1}N_{\ell k}$, our assumption that $\alpha^{-1}M_{i j} > -1$ implies that $\beta^{-1}N_{\ell k} < -1$, and this holds for all $\ell \sim k$. However, this implies that $\beta^{-1}N$ is the Gram matrix of a vector $t$-coloring in which every edge is slack. This implies that $\chiv(H) < t$, a contradiction. Therefore, we must have that $\alpha^{-1} M_{ij} \le -1$ for all $i \sim j$, and thus $\alpha^{-1} M$ is the Gram matrix of a vector $t$-coloring of $G$. Symmetrically, $\beta^{-1}N$ is the Gram matrix of a vector $t$-coloring of $H$. Finally, we have that $M \otimes J + J \otimes N$ is a convex combination of $\alpha^{-1} M \otimes J$ and $J \otimes \beta^{-1}N$ with coefficients $\alpha$ and $\beta$ respectively.\qeds

The above lemma and those following will be used to prove Result~\ref{res:A}, which is Theorem~\ref{thm:A} below. First we prove the lower bound on $\rk(G \times H)$ in the $\chiv(G) = \chiv(H)$ case.

\begin{lemma}\label{lem:rkbound}
Let $G$ and $H$ be graphs with $\chiv(G) = \chiv(H)$. Then any convex combination of optimal vector colorings induced by $G$ and $H$ has rank at most $\rk(G) + \rk(H)$, and equality can be attained. Thus, $\rk(G \times H) \ge \rk(G) + \rk(H)$.
\end{lemma}
\proof
If either factor is empty then all of the ranks are zero and the lemma holds trivially. So we may assume that $G$ and $H$ are non-empty.

Consider an arbitrary convex combination of Gram matrices of optimal vector colorings induced by $G$ and $H$. This has the form
\[\sum_{i = 1}^m \alpha_i M_i \otimes J + \sum_{j=1}^n \beta_j J \otimes N_j,\]
where $M_1, \ldots, M_m$ and $N_1, \ldots, N_n$ are Gram matrices of optimal vector colorings of $G$ and $H$ respectively, and $\alpha_i , \beta_j \ge 0$ for all $i,j$, and $\sum_i \alpha_i + \sum_j \beta_j = 1$. We can rewrite the convex combination above as $\alpha M \otimes J + \beta J \otimes N$ where
\[M = \frac{1}{\sum_i \alpha_i}\left(\sum_i \alpha_i M_i\right), \quad N = \frac{1}{\sum_j \beta_j}\left(\sum_j \beta_j N_j\right),\]
and $\alpha  = \sum_i \alpha_i$ and $\beta = \sum_j \beta_j$. In other words, we can always reduce to the case of a convex combination of a single vector coloring induced by $G$ and a single vector coloring induced by $H$.

In this case, since $\rk(J) = 1$, it is easy to see that
\[\rk(\alpha M \otimes J + \beta J \otimes N) \le \rk(M) + \rk(N) \le \rk(G) + \rk(H).\]
Thus we have proven the first inequality, and it only remains to show that equality can be attained.

Let $\p$ and $\q$ be max-rank vector colorings of $G$ and $H$ respectively. Consider the optimal vector coloring $\w$ of $G \times H$ which is given by $w_{i\ell} = (1/\sqrt{2})(p_i \oplus q_\ell)$, i.e., $\w$ is a direct sum of vector colorings of $G \times H$ induced by $G$ and $H$.

Since $G$ and $H$ are nonempty, we have that $\chiv(G) = \chiv(H) \ge 2$, and then by Lemma~\ref{lem:skcomps} both $G^{\sk}$ and $H^{\sk}$ are nonempty. Thus both $G$ and $H$ contain some neighborly vertices. Therefore, by Lemma~\ref{lem:cone2conv} there exist nonnegative coefficients $\delta_i$ for $i \in V(G)$ such that
\begin{equation*}\label{eq:scale}
\sum_i \delta_i p_i = \zero \quad \& \quad \sum_i \delta_i = 1.
\end{equation*}
Thus we have that
\begin{equation*}\label{eq:deltai}
\sum_i \delta_i (p_i \oplus q_\ell) = \zero \oplus q_\ell.
\end{equation*}
This gives us $\zero \oplus q_\ell \in \spn(\w)$ for all $\ell \in V(H)$. Similarly, $p_i \oplus \zero \in \spn(\w)$ for all $i \in V(G)$. Using these vectors we can obtain any vector in $\spn(\p) \oplus \spn(\q)$ and therefore $\rk(G \times H) \ge \dim \spn(\w) = \rk(G) + \rk(H)$.\qeds

In the proof of Theorem~\ref{thm:A}, we will need to make use of Lemma~\ref{lem:conespan}. In order to be able to do this, we will need the following:

\begin{lemma}\label{lem:rk2nbrly}
Let $G$ and $H$ be non-empty graphs with $\chiv(G) = \chiv(H)$. If $\rk(G \times H) = \rk(G) + \rk(H)$, then all vertices of $G$ and $H$ are neighborly.
\end{lemma}
\proof
First, we will show that $(G \times H)^{\sk}$ is a (spanning) subgraph of $G^{\sk} \times H^{\sk}$. By Lemma~\ref{lem:sksandwich} there exist optimal vector colorings $\p$ and $\q$ of $G$ and $H$ respectively such that $G^\p = G^{\sk}$ and $H^\q = H^{\sk}$. Letting $\w$ be the optimal vector coloring of $G \times H$ defined as $w_{i\ell} = (1/\sqrt{2})(p_i \oplus q_\ell)$, it is easy to see that $(G \times H)^\w = G^\p \times H^\q = G^{\sk} \times H^{\sk}$, and thus $(G \times H)^{\sk}$ is a subgraph of $G^{\sk} \times H^{\sk}$ (this actually does not depend on $\rk(G \times H) = \rk(G) + \rk(H)$).

Next we will show that $\rk(G \times H) \ge \rk(G^{\sk} \times H^{\sk})$. First, let $t = \chiv(G) = \chiv(H)$ and note that $\chiv(G^{\sk} \times H^{\sk}) = t = \chiv((G \times H)^{\sk})$ by Lemma~\ref{lem:skcomps} and Theorem~\ref{thm:hedetniemi}. Now, since $(G \times H)^{\sk}$ is a spanning subgraph of $G^{\sk} \times H^{\sk}$, every optimal vector coloring of the latter is an optimal vector coloring of the former. Therefore $\rk((G \times H)^{\sk}) \ge \rk(G^{\sk} \times H^{\sk})$, and so $\rk(G \times H) \ge \rk(G^{\sk} \times H^{\sk})$ by Lemma~\ref{lem:skcomps}.

Now suppose that some vertex $i^*$ of $G$ is not neighborly and therefore is isolated in $G^{\sk}$. Let $G'$ be the graph obtained from $G^{\sk}$ by removing $i^*$. Thus $G^{\sk} \cong G' \cup K_1$ and $\rk(G^{\sk}) = \rk(G') + 1$. Using this we can rewrite $G^{\sk} \times H^{\sk}$ as $(G' \times H^{\sk}) \cup (K_1 \times H^{\sk})$. Note that $K_1 \times H^{\sk}$ is simply $|V(H)|$ isolated vertices. Every isolated vertex adds exactly one to the rank of a graph and so $\rk(G^{\sk} \times H^{\sk}) = \rk(G' \times H^{\sk}) + |V(H)| \ge \rk(G' \times H^{\sk}) + 2$ since $H$ is nonempty. However, $\rk(G') = \rk(G) - 1$ and so $\rk(G' \times H^{\sk}) \ge \rk(G) - 1 + \rk(H)$ by Lemma~\ref{lem:rkbound}. Combining all of this we have that
\[\rk(G \times H) \ge \rk(G^{\sk} \times H^{\sk}) \ge \rk(G' \times H^{\sk}) + 2 \ge \rk(G) + \rk(H) + 1,\]
a contradiction. Therefore every vertex of $G$, and similarly of $H$, is neighborly.\qeds

We remark that it is not much more difficult to show that $G^{\sk}$ and $H^{\sk}$ must be connected, but we do not need this here and so we save it for Corollary~\ref{cor:skconnected}. The last lemma we will need shows that we can use dependencies among the vectors in vector colorings of $G$ and $H$ to construct dependencies among the vectors in a vector coloring of $G \times H$.

\begin{lemma}\label{lem:dependencies}
Let $G$ and $H$ be graphs with $\chiv(G) = \chiv(H)$, and let $\p$ and $\q$ be optimal vector colorings of $G$ and $H$ respectively. Let $\w$ be the optimal vector coloring of $G \times H$ given by $w_{i\ell} = (1/\sqrt{2}) (p_i \oplus q_\ell)$ for $i \in V(G)$, $\ell \in V(H)$. If $i \to_\p j$ and $\ell \to_\q k$, then $(i,\ell) \to_\w (j,k)$.
\end{lemma}
\proof
Let $t = \chiv(G) = \chiv(H)$. If $i \to_\p j$ and $\ell \to_\q k$, then there are nonnegative coefficients $\alpha_{j'}$ for $j' \in N^\p(i)$ and $\beta_{k'}$ for $k' \in N^\q(\ell)$ such that
\[p_i + \sum_{j' \in N^\p(i)} \alpha_{j'} p_{j'} = \zero \quad \text{and} \quad q_\ell + \sum_{k' \in N^\q(\ell)} \beta_{k'} q_{k'} = \zero,\]
and $\alpha_{j}, \beta_{k} > 0$. Moreover, by Lemma~\ref{lem:cone2conv} we have that
\[\sum_{j' \in N^\p(i)} \alpha_{j'} = t-1 = \sum_{k' \in N^\q(\ell)} \beta_{k'}.\]
Now, we have that
\begin{align*}
\zero =& \left(p_i \oplus \zero + \sum_{j' \in N^\p(i)} (\alpha_{j'} p_{j'}) \oplus \zero\right) + \left(\zero \oplus q_\ell + \sum_{k' \in N^\q(\ell)} \zero \oplus (\beta_{k'} q_{k'})\right) \\
=& \ p_i \oplus q_\ell + \left(\sum_{j' \in N^\p(i)} \frac{1}{t-1} \sum_{k' \in N^\q(\ell)} (\beta_{k'} \alpha_{j'} p_{j'}) \oplus \zero\right) \\
&+ \left(\sum_{k' \in N^\q(\ell)} \frac{1}{t-1} \sum_{j' \in N^\p(i)} \zero \oplus (\alpha_{j'}\beta_{k'} q_{k'})\right) \\
=& \ p_i \oplus q_\ell + \frac{1}{t-1} \left(\sum_{j' \in N^\p(i), k' \in N^\q(\ell)} (\beta_{k'}\alpha_{j'} p_{j'}) \oplus (\alpha_{j'}\beta_{k'} q_{k'}) \right) \\
=& \ p_i \oplus q_\ell + \frac{1}{t-1} \sum_{j' \in N^\p(i), k' \in N^\q(\ell)} \alpha_{j'}\beta_{k'} (p_{j'} \oplus q_{k'}).
\end{align*}
Since $w_{i\ell} = (1/\sqrt{2}) (p_i \oplus q_\ell)$ and $N^\w(i,\ell) = N^\p(i) \times N^\q(\ell)$, this implies that
\[w_{i\ell} + \frac{1}{t-1}\sum_{(j',k') \in N^\w(i,\ell)} \alpha_{j'}\beta_{k'} w_{j'k'} = 0.\]
Moreover, we have that the coefficient of $w_{jk}$ is $\alpha_{j}\beta_{k}/(t-1) > 0$. Therefore, $(i,\ell) \to_\w (j,k)$ as desired.\qeds

Finally, we can prove Result~\ref{res:A}:

\begin{theorem}\label{thm:A}
Let $G$ and $H$ be graphs such that $\chiv(G) = \chiv(H)$. Then $\rk(G \times H) = \rk(G) + \rk(H)$ if and only if every optimal vector coloring of $G \times H$ is a convex combination of vector colorings induced by $G$ and $H$.
\end{theorem}
\proof
By Lemma~\ref{lem:rkbound}, the maximum rank attained by any optimal vector of $G \times H$ that is a convex combination of vector colorings induced by $G$ and $H$ is $\rk(G) + \rk(H)$. Therefore, if $\rk(G \times H) > \rk(G) + \rk(H)$, then some optimal vector coloring of $G \times H$ does not have this form. So we have proven one direction of the claim.

Conversely, suppose that $\rk(G \times H) = \rk(G) + \rk(H)$. Note that if either factor is empty then both factors and the product are empty and the claim follows trivially. Thus we may assume that both $G$ and $H$ are non-empty.

Let $\p$ and $\q$ be optimal vector colorings of $G$ and $H$ that span $\mathbb{R}^{\rk(G)}$ and $\mathbb{R}^{\rk(H)}$ respectively. Define $\w$ to be the optimal vector coloring of $G \times H$ given by $w_{i\ell} = (1/\sqrt{2})(p_i \oplus q_\ell)$, and let $W$ be the matrix whose rows are the vectors in $\w$. Note that, as in the proof of Lemma~\ref{lem:rkbound}, the vector coloring $\w$ has rank $\rk(G) + \rk(H)$, i.e., it spans the space $\mathbb{R}^{\rk(G)+\rk(H)}$ which it is contained in. Therefore, $\w$ is a max-rank vector coloring of $G \times H$. Thus, by Theorem~\ref{thm:main} we have that the Gram matrix of any optimal vector coloring of $G \times H$ is equal to $W(I+R)W^T$ for some symmetric matrix $R$ satisfying
\begin{equation}\label{eq:Rw}
\begin{aligned}
&w_{i\ell}^T R w_{i\ell} = 0 \text{ for all } i,\ell; \\
&w_{i\ell}^T R w_{jk} \le -1 - w_{i\ell}^T w_{jk} \text{ for } (i,\ell) \sim (j,k); \\
&I+R \succeq 0,
\end{aligned}
\end{equation}
We will show that for any such $R$, the matrix $W(I+R)W^T$ is a convex combination of vector colorings induced by $G$ and $H$, thus proving the theorem.

Suppose we have a symmetric matrix $R$ satisfying~\eqref{eq:Rw}, and partition it into block form according to the partition of the $\w$ vectors with respect to $\p$ and $\q$:
\[R = \begin{pmatrix}
R_1 & F \\
F^T & R_2
\end{pmatrix}.\]
We will begin by showing that $F = 0$. The first step is to show that $i \to_\p j$ and $\ell \to_\q k$ implies that $(p_i - p_j)^T F (q_\ell - q_k)\le 0.$ Thus, suppose that $i \to_\p j$ and $\ell \to_\q k$. In particular, we have that $i \sim j$ and $\ell \sim k$. By Lemma~\ref{lem:dependencies}, we have that $(i,\ell) \to_\w (j,k)$. Since $\w$ is a max-rank vector coloring of $G \times H$, Lemma~\ref{lem:nbrlymr} implies that $(i,\ell) \to (j,k)$ and thus $(i,\ell) \sim_{\sk} (j,k)$ (recall Remark~\ref{rem:arrow2sk}). As noted in Remark~\ref{rem:tightsk}, this implies that $w_{i\ell}^T R w_{jk} = 0$. Therefore,
\begin{equation}\label{eq:Feq}
0 = 2w_{i\ell}^T R w_{jk} = \left(p_i^T R_1 p_j + q_\ell^T R_2 q_k\right) + \left(p_i^T F q_k + p_j^T F q_\ell \right).
\end{equation}
Next we consider the vertices $(i,k)$ and $(j,\ell)$. In this case we do not know that $(i,k) \to (j,\ell)$ since we do not know whether $\ell \to k$ implies $k \to \ell$. However, since $(i,\ell) \to_\w (j,k)$,  we have that  $(i,\ell) \sim_\w (j,k)$, and thus
\[w_{ik}^Tw_{j\ell} = \frac{1}{2} (p_i^Tp_j + q_k^Tq_\ell) = \frac{1}{2} (p_i^Tp_j + q_\ell^Tq_k) = w_{i\ell}^Tw_{jk} = -1.\]
Therefore, since $(i,k) \sim (j,\ell)$, we have
\begin{equation}\label{eq:Fineq}
0 \ge 2w_{ik}^T R w_{j\ell} = \left(p_i^T R_1 p_j + q_\ell^T R_2 q_k\right) + \left(p_i^T F q_\ell + p_j^T F q_k \right),
\end{equation}
where we have used the fact that $R_2$ is symmetric. Subtracting (\ref{eq:Feq}) from (\ref{eq:Fineq}), we obtain
\[0 \ge p_i^T F q_\ell + p_j^T F q_k - p_i^T F q_k - p_j^T F q_\ell = (p_i - p_j)^T F (q_\ell - q_k).\]
This holds for all $i \to_\p j$ and $\ell \to_\q k$. Therefore we have that $p^T F q \le 0$ for all $p \in \cone(\{p_i-p_j : i \to_\p j\})$ and $q \in \cone(\{q_\ell - q_k : \ell \to_\q k\})$. But since we assumed that $\rk(G \times H) = \rk(G) + \rk(H)$, Lemma~\ref{lem:rk2nbrly} tells us that every vertex of $G$ and $H$ is neighborly. Therefore, by Lemma~\ref{lem:conespan}, we have that these two cones are equal to $\spn(\p)$ and $\spn(\q)$ respectively. This implies that $F = 0$ as desired.

Now, note that the matrix $W^T$ can be written as
\[W^T = \frac{1}{\sqrt{2}}\begin{pmatrix}
P^T \otimes \one^T \\
\one^T \otimes Q^T
\end{pmatrix}\]
where $P$ and $Q$ are the matrices whose rows are the vectors in $\p$ and $\q$ respectively. Using the fact that $I + R_1 = (I+R_1) \otimes I_1$ and $I+R_2 = I_1 \otimes (I+R_2)$ where $I_1$ is the $1 \times 1$ identity matrix, we see that
\begin{align*}
W(I+R)W^T &= \frac{1}{2}\begin{pmatrix}P \otimes \one \ \ \one \otimes Q \end{pmatrix} \begin{pmatrix}I+R_1 & 0 \\ 0 & I+R_2\end{pmatrix} \begin{pmatrix}P^T \otimes \one^T  \\ \one^T \otimes Q^T\end{pmatrix} \\
&= \frac{1}{2} \Big(\left[P(I+R_1)P^T \otimes \one \one^T\right] + \left[\one \one^T \otimes Q(I+R_2)Q^T\right]\Big) \\
&= \frac{1}{2} \Big(\left[P(I+R_1)P^T \otimes J\right] + \left[J \otimes Q(I+R_2)Q^T\right]\Big).
\end{align*}
Since $I+R \succeq 0$, we have that $I+R_i \succeq 0$ for $i = 1,2$. Thus the matrices $P(I+R_1)P^T$ and $Q(I+R_2)Q^T$ are positive semidefinite. Therefore, by Lemma~\ref{lem:convcombvc}, the matrix $W(I+R)W^T$ is the Gram matrix of a convex combination of vector colorings induced by $G$ and $H$.\qeds

Unfortunately, the hypothesis $\rk(G \times H) = \rk(G) + \rk(H)$ of Theorem~\ref{thm:A} depends on both the product and the factors. It would be preferable to obtain a similar characterization where the hypotheses depended only on the factors individually, but we were not able to prove such a result. In fact, it may not be possible. If there exist graphs $G_1,G_2$ and $H_1,H_2$ with all the same vector chromatic number and such that $\rk(G_1 \times H_2) = \rk(G_1) + \rk(H_2)$ and $\rk(G_2 \times H_1) = \rk(G_2) + \rk(H_1)$, but $\rk(G_1 \times H_1) > \rk(G_1) + \rk(H_1)$, then any characterization must take into account some property of the pair of factors, and not just properties of the factors individually. However, we were not able to find such graphs. We are able to prove an analog of Corollary~\ref{cor:B} which provides a sufficient condition based only on the individual factors:

\begin{corollary}\label{cor:A}
Let $G$ and $H$ be graphs such that $\chiv(G) = \chiv(H)$. If both $G$ and $H$ admit connected strictly complementary dual solutions, then every optimal vector coloring of $G \times H$ is a convex combination of vector colorings induced by the factors.
\end{corollary}
\proof

Let $t = \chiv(G) = \chiv(H)$, and let $B_G$ and $B_H$ be connected strictly complementary dual solutions for $G$ and $H$ respectively. Recall that this means that $\text{corank}(B_G) = \rk(G)$ and $\text{corank}(B_H) = \rk(H)$. Since $B_G$ and $B_H$ are connected and positive semidefinite, they have strictly positive diagonal and so we can apply Lemma~\ref{lem:dualconvert2} to obtain optimal solutions $A_G$ and $A_H$ to~(\ref{eq:Adual}) for $G$ and $H$. These solutions have the property that $G(A_G) = G(B_G)$ and the multiplicity of $-1$ as an eigenvalue of $A_G$ is equal to $\text{corank}(B_G) = \rk(G)$, and similarly for $A_H$. This implies that both $A_G$ and $A_H$ are connected and therefore they have strictly positive maximum eigenvectors and their maximum eigenvalues have multiplicity one. Since they are optimal solutions to~(\ref{eq:Adual}), the maximum eigenvalue of both $A_G$ and $A_H$ is $\lambda = t-1 \ge 1$. We will consider the matrix $A = \frac{1}{\lambda} A_G \otimes A_H$. Recall from the proof of Theorem~\ref{thm:hedetniemi} that $A$ is an optimal solution to~(\ref{eq:Adual}) for $G \times H$ and $A$ has maximum eigenvalue $\lambda$. We can construct an eigenvector for this eigenvalue by taking the Kronecker product of the strictly positive maximum eigenvectors of $A_G$ and $A_H$. Thus $A$ also has a strictly positive maximum eigenvector. Finally, since the maximum eigenvalues of both $A_G$ and $A_H$ are simple, the multiplicity of $-1$ as an eigenvalue of $A$ is equal to the sum of the multiplicities of $-1$ as an eigenvalue of $A_G$ and $A_H$, and this is $\rk(G) + \rk(H)$.

Now we can apply Lemma~\ref{lem:dualconvert2} again to obtain an optimal dual solution $B$ for $\chiv(G \times H)$ with corank $\rk(G) + \rk(H)$. This implies that $\rk(G \times H) \le \rk(G) + \rk(H)$, but of course we already have the other inequality. Therefore, by Theorem~\ref{thm:A}, we have proven the corollary.\qeds

We immediately obtain to the following corollary for the case where $G$ and $H$ are uniquely vector colorable:

\begin{corollary}\label{cor:Auniq}
Let $G$ and $H$ be graphs such that $\chiv(G) = \chiv(H)$, and which both admit connected strictly complementary dual solutions. Suppose that $G$ and $H$ have unique vector colorings $\p$ and $\q$ respectively. Then the only vector colorings of $G \times H$ are, up to isometry, the direct sums of the vector colorings induced by $\p$ and $\q$.
\end{corollary}

Note that we needed to assume that the strictly complementary dual solutions for $G$ and $H$ in Corollary~\ref{cor:A} above were connected, not just that they had strictly positive diagonal as in Corollary~\ref{cor:B}. In fact it is not difficult to show that if one tries to construct the optimal dual solution $B$ in the proof above starting from $B_G$ and $B_H$ that are not both connected, then the corank of $B$ will be greater than $\rk(G) +\rk(H)$, and so this will not suffice to obtain the conclusion. Of course, this does not prove that $G$ and $H$ having connected strictly complementary dual solutions is necessary, but as in Section~\ref{subsec:Btheorem}, we can show that an analog of this property in terms of skeletons is necessary:

\begin{corollary}\label{cor:skconnected}
Let $G$ and $H$ be nonempty graphs such that $\chiv(G) = \chiv(H)$. If every optimal vector coloring of $G \times H$ is a convex combination of vector colorings induced by the factors, then $G^{\sk}$ and $H^{\sk}$ are connected.
\end{corollary}
\proof
Let $t = \chiv(G) = \chiv(H)$. Since every optimal vector coloring of $G \times H$ is a convex combination of vector colorings induced by the factors, we have that $\rk(G \times H) = \rk(G) + \rk(H)$. Recall from Lemma~\ref{lem:rk2nbrly} that this implies that every vertex of $G$ and $H$ is neighborly. By Lemma~\ref{lem:nbrlysk} this means that there are no isolated vertices in $G^{\sk}$ or $H^{\sk}$. Let $G_1, \ldots, G_r$ and $H_1, \ldots, H_s$ be the connected components of $G^{\sk}$ and $H^{\sk}$ respectively. By the above we have that none of these components are empty and by Lemma~\ref{lem:skcomps} we have that they all have vector chromatic number equal to $t$. Therefore, by Lemma~\ref{lem:rkbound}, we have that $\rk(G_i \times H_j) \ge \rk(G_i) + \rk(H_j)$ for all $i \in [r]$, $j \in [s]$. Recall also from the proof of Lemma~\ref{lem:rk2nbrly} that $\rk(G \times H) \ge \rk(G^{\sk} \times H^{\sk})$ (in fact we have equality in this case since it is easy to see that $(G \times H)^{\sk} = G^{\sk} \times H^{\sk}$ if every optimal vector coloring of $G \times H$ is a convex combination of vector colorings induced by the factors). Using the fact that $\rk(G) = \sum_i \rk(G_i)$ and $\rk(H) = \sum_j \rk(H_j)$ from Lemma~\ref{lem:skcomps}, we see that
\begin{align*}
\rk(G \times H) &\ge \rk(G^{\sk} \times H^{\sk}) \\
&= \sum_{i=1}^r \sum_{j=1}^s \rk(G_i \times H_j) \\
&\ge \sum_{i=1}^r \sum_{j=1}^s \left(\rk(G_i) + \rk(H_j)\right) \\
&= s\sum_{i=1}^r \rk(G_i) + r\sum_{j=1}^s \rk(H_j) \\
&= s\rk(G) + r\rk(H).
\end{align*}
Obviously, if either $r$ or $s$ is greater than 1, then the last expression is strictly greater than $\rk(G) + \rk(H)$. Thus, $G^{\sk}$ and $H^{\sk}$ must be connected.\qeds

Combining Corollaries~\ref{cor:A} and~\ref{cor:skconnected}, we obtain the following:

\begin{theorem}\label{thm:Asandwich}
Let $G$ and $H$ be nonempty graphs with $\chiv(G) = \chiv(H)$. Then the following statements satisfy $(1) \Rightarrow (2) \Rightarrow (3)$:
\begin{enumerate}
\item Both $G$ and $H$ admit connected strictly complementary dual solutions.
\item Every optimal vector coloring of $G \times H$ is a convex combination of optimal vector colorings induced by the factors.
\item Both $G^{\sk}$ and $H^{\sk}$ are connected.
\end{enumerate}
\end{theorem}

Whether or not any of the above implications can be reversed is an interesting open question.

\subsection{1-Walk-Regular Graphs}\label{subsec:1walkreg}

A graph $G$ is 1-walk-regular if for all $k \in \mathbb{N}$,
\begin{enumerate}
\item the number of walks of length $k$ starting and ending at a vertex of $G$ is independent of the vertex;
\item the number of walks of length $k$ starting at one end of an edge and ending at the other is independent of the edge.
\end{enumerate}
This definition can be written algebraically as follows: A graph with adjacency matrix $A$ is 1-walk-regular if there exist $a_k, b_k \in \mathbb{N}$ for all $k \in \mathbb{N}$ such that
\begin{enumerate}
\item $A^k \circ I = a_k I$;
\item $A^k \circ A = b_k A$.
\end{enumerate}
Note that any 1-walk-regular graph must be regular. Also, any graph which is vertex- and edge-transitive is easily seen to be 1-walk-regular. Other classes of 1-walk-regular graphs include distance regular graphs and, more generally, graphs which are a single class in an association scheme.

Suppose $G$ is 1-walk-regular with adjacency matrix $A$ which has least eigenvalue $\tau$. Let $E_\tau$ be the projector onto the $\tau$-eigenspace of $A$. In~\cite{UVC1}, we showed that $E_\tau$ and $A - \tau I$ are, up to positive scalars, optimal primal and dual solutions for $\chiv(G)$ respectively. Since $\ker(A - \tau I) = \im(E_\tau)$, these form a strictly complementary pair for $G$. The matrix $A - \tau I$ clearly has strictly positive diagonal, and is connected if and only if the graph $G$ is connected. Because of this, the statements $(A')$ and $(B')$ can be proven for connected 1-walk-regular graphs, and in fact follow from Corollaries~\ref{cor:Auniq} and~\ref{cor:B2}. Specifically, we have the following:

\begin{theorem}
If $G$ is a 1-walk-regular graph, then the following hold:
\begin{enumerate}
\item If $H$ is a connected graph and $\chiv(G) < \chiv(H)$, then every optimal vector coloring of $G \times H$ is induced by $G$.
\item If, additionally, $G$ is uniquely vector colorable, then $G \times H$ is uniquely vector colorable.
\end{enumerate}
If both $G$ and $H$ are connected 1-walk-regular graphs, then the following hold:
\begin{enumerate}
\setcounter{enumi}{2}
\item Every optimal vector coloring of $G \times H$ is a convex combination of vector colorings induced by $G$ and $H$.
\item If, additionally, $G$ and $H$ have unique vector colorings $\p$ and $\q$ respectively, then the only optimal vector colorings of $G \times H$ are the convex combinations of the vector colorings induced by $\p$ and $\q$.
\end{enumerate}
\end{theorem}

Statements $(1)$ and $(3)$ in the theorem above follow from Corollaries~\ref{cor:A} and~\ref{cor:B} respectively, and they imply statements $(2)$ and $(4)$ respectively.

\subsection{Implications}\label{subsec:implications}

As we mentioned in Section~\ref{subsec:motivations}, Duffus et.~al.~showed that $(A_n) \Rightarrow (B_n) \Rightarrow (C_{n+1})$ for all positive $n$. If we similarly parameterize the statements $(A')$, $(B')$, and $(C')$, do the same implications hold? Firstly, we cannot parameterize these statements using only integers, since $\chiv$ can take on non-integer, and even irrational, values. Since, other than for the empty graph, the vector chromatic number is always at least 2, we parameterize the statements using real numbers $t \ge 2$:
\begin{itemize}
\item $\boldsymbol{(A'_t):}$ For all uniquely vector $t$-colorable graphs $G$ and $H$, each vector $t$-coloring of $G \times H$ is a convex combination of the vector $t$-colorings induced by $G$ and $H$.
\item $\boldsymbol{(B'_t):}$ For all uniquely vector $t$-colorable graphs $G$ and connected graphs $H$ with $\chiv(H) > t$, the graph $G \times H$ is uniquely vector $t$-colorable.
\item $\boldsymbol{(C'_t):}$ For all graphs $G$ and $H$ with $\min\{\chiv(G),\chiv(H)\} = t$, we have $\chiv(G \times H) = t.$
\end{itemize}
Note the difference between $(C'_t)$ and $(C_n)$. In $(C_n)$ it was assumed that both graphs have chromatic number $n$, instead of just the minimum being $n$. However this is just as general since any graph with chromatic number at least $n$ has a subgraph of chromatic number exactly $n$. But the same is not true of graphs with vector chromatic number at least $t$, so we need the more general statement here if we want $(C')$ to be equivalent to $(C'_t)$ being true for all $t$.

Another problem that arises is that the relevant relationship between $n$ and $n+1$ for proving $(B_n) \Rightarrow (C_{n+1})$ is that $n+1$ is the smallest value achievable by $\chi$ that is greater than $n$. It is known that the vector chromatic number of the Kneser graph $K_{n:r}$ is $n/r$~\cite{Lovasz}, so there are no two ``consecutive" values for $\chiv$. Of course, since we have already proven $(C')$, the implication $(B'_t) \Rightarrow (C'_{t+1})$ technically does hold since $(C'_{t+1})$ is always true. Since this implication is trivial, we will focus on the implication $(A'_t) \Rightarrow (B'_t)$.

A key ingredient in the proof of $(A_n) \Rightarrow (B_n)$ is the result of Greenwell and Lov\'{a}sz~\cite{greenwell} which says that if a graph $H$ is connected and $\chi(H) > n$, then $H \times K_n$ is uniquely $n$-colorable. In order to adapt their proof, we would need an analog of this result, which would include an analog of $K_n$ for every $t$ which is the vector chromatic number of some graph. Since it is not even known what real numbers can be obtained as the vector chromatic number of a graph, this seems difficult to do. Therefore, we define the following statement for all $t$ that is the vector chromatic number of some graph:

\vspace{.1cm}
$\boldsymbol{(D'_t):}$ There exists a graph $G_t$ such that if $H$ is connected and $\chiv(H) > t$, then $H \times G_t$ is uniquely vector $t$-colorable.
\vspace{.1cm}

Note that for all rational $t = n/r$, letting $G_t = K_{n:r}$ works for the above statement. Also, by Corollary~\ref{cor:B2}, to prove $(D'_t)$ it suffices to find a uniquely vector $t$-colorable graph for which there exists a strictly complementary dual solution with strictly positive diagonal. The following theorem shows that $(A'_t) \Rightarrow (B'_t)$ if we assume $(D'_t)$.

\begin{theorem}\label{thm:implication}
For all $t \in \mathbb{R}$, we have that $(A'_t)$ \& $(D'_t)$ implies $(B'_t)$.
\end{theorem}
\proof
We follow the proof of Duffus et.~al.~for chromatic number. Suppose that $(A'_t)$ and $(D'_t)$ hold, that $G$ is uniquely vector $t$-colorable, and that $H$ is connected with $\chiv(H) > t$. Let $G_t$ be the graph guaranteed by $(D'_t)$. If $G \times H$ is uniquely vector $t$-colorable, then we are done. Otherwise $G \times H$ has at least two vector $t$-colorings. One of these is the vector coloring induced by the unique vector coloring of $G$. This vector coloring is independent of $H$, meaning that the $(i,\ell)(j,k)$-entry of its Gram matrix is determined by $i$ and $j$ alone. The other vector coloring of $G \times H$ cannot have this property, since this would give a distinct vector coloring of $G$.

Now consider the graph $G \times H \times G_t$. Written as $G \times (H \times G_t)$, the statements $(A'_t)$ and $(D'_t)$ imply that the only vector colorings this graph has are the convex combinations of the vector colorings induced by $G$ and $H \times G_t$. Note that both of these induced colorings are independent of $H$ in the above described manner, since the unique vector coloring of $H \times G_t$ is the one induced by $G_t$. This means that their convex combinations are also independent of $H$. However, written as $(G \times H) \times G_t$, this graph also has a vector coloring induced by the vector coloring of $G \times H$ that is not independent of $H$, a contradiction.\qeds

\section{Discussion}\label{sec:discuss}

The main results of this paper were a vector coloring analog of Hedetniemi's Conjecture (Theorem~\ref{thm:hedetniemi}), and Theorems~\ref{thm:B} and~\ref{thm:A} showing that the optimal vector colorings of the categorical product of two graphs can be described in terms of the optimal vector colorings of the factor(s) with the minimum vector chromatic number if and only if certain conditions hold. Though the conditions in Theorems~\ref{thm:B} and~\ref{thm:A} are necessary and sufficient, the conditions concern both the product graph and the factors, rather than just the factors. It would be a significant improvement to obtain similar results but whose conditions depend only on the factors. Possible such conditions are given by the first and last items of Theorems~\ref{thm:Bsandwich} and~\ref{thm:Asandwich} respectively.

Theorems~\ref{thm:B} and~\ref{thm:A} do not require unique vector colorability of the factor(s), but they also do not imply statements $(B')$ and $(A')$ respectively. Thus it is an open question as to whether either of these statements are true. In order to prove these it would be useful to find some necessary conditions for a graph to be uniquely vector colorable. Following ideas along the lines of those in Section~\ref{subsec:specialprops}, it is not hard to see that a nonempty uniquely vector colorable graph must have a connected skeleton. This is good since Corollary~\ref{cor:skconnected} tells us that this is necessary for the conclusion of statement $(A')$ to hold. Can we show that any uniquely vector colorable graph must admit a strictly complementary dual solution? We know of no counterexample. In fact we do not know a single graph which we know does not admit a strictly complementary dual solution. Even proving this would not be quite enough, since we need the strictly complementary dual to have strictly positive diagonal (in the case of statement $(B')$) or be connected (in the case of statement $(A')$) in order to apply Corollary~\ref{cor:B} or~\ref{cor:A}.

There are many questions one could ask regarding skeletons. Lemma~\ref{lem:sksandwich} points out that if $B$ is an optimal dual solution for a graph $G$, then $G(B)$ is a subgraph of $G^{\sk}$. Is it possible that there always exists an optimal dual solution $B$ such that $G(B) = G^{\sk}$? Perhaps the notion of neighborliness can help here. Using ideas similar to those in the proof of Lemma~\ref{lem:nbrlymr}, it is not difficult to show that for any neighborly vertex $i$, there exists a fixed convex combination with support $D[i]$ that witnesses the neighborliness of $i$ in every optimal vector coloring. Can these convex combinations be used to construct an optimal dual solution $B$ such that $B_{ij} > 0$ if $j \in D[i]$? This would not quite be enough, since we do not know that one of $j \in D(i)$ or $i \in D(j)$ holds for every edge $ij$ in $G^{\sk}$, but it would be a big step in the right direction.

In the proof of Lemma~\ref{lem:rk2nbrly}, we showed that $(G \times H)^{\sk}$ is a subgraph of $G^{\sk} \times H^{\sk}$ whenever $\chiv(G) = \chiv(H)$. If every optimal vector coloring of $G \times H$ is a convex combination of vector colorings induced by $G$ and $H$, then it is easy to see that $(G \times H)^{\sk} = G^{\sk} \times H^{\sk}$. This feels analogous to the situation with the inequality $\rk(G \times H) \ge \rk(G) + \rk(H)$. This begs the question of whether an analog of Theorem~\ref{thm:A} holds for skeletons, i.e., whether $(G \times H)^{\sk} = G^{\sk} \times H^{\sk}$ if and only if every optimal vector coloring of $G \times H$ is a convex combination of vector colorings induced by $G$ and $H$. Similarly, in the $\chiv(G) < \chiv(H)$ case, does $(G \times H)^{\sk} = G^{\sk} \times H$ if and only if every optimal vector coloring of $G \times H$ is induced by $G$?

As mentioned above, we do not know a single graph which we know has no strictly complementary dual solution. Thus the existence of such a graph is an important open question. One approach to this would be to consider graphs that lie in (symmetric) association schemes. In this case, the primal and dual semidefinite programs in \eqref{eq:primal} and (\ref{eq:dual}) always have optimal solutions that lie in the association scheme, and thus they can be replaced by linear programs. Moreover, given any optimal solution to the primal or dual, one can project this solution down into the association scheme while preserving its optimality and without decreasing its rank (this is nontrivial). Therefore, such a graph has a strictly complementary pair of primal and dual solutions if and only if it has such a pair lying in the association scheme. Any element of the association scheme is a linear combination of the idempotents of the scheme. Therefore, in order to show that such a graph does not satisfy strict complementarity, it suffices to find an idempotent which must have a coefficient of zero in any linear combination giving an optimal primal or dual solution. For a fixed idempotent this can be checked by solving a linear program which maximizes the coefficient of the idempotent in any such a linear combination for the primal/dual (the vector chromatic number must be computed first in order to construct this LP). If the max value is zero for both the primal and the dual for a given idempotent, then the graph does not satisfy strict complementarity. Note that we can restrict to graphs that are not a single class in an association scheme for this approach since otherwise the graph will be 1-walk-regular which implies it will satisfy strict complementarity.

Though we have not written down detailed proofs, we believe that the \emph{strict} vector coloring analogs of Results~\ref{res:B} and~\ref{res:A} hold. In fact, we believe that the proofs of these results should be much simpler than the results here, because for strict vector colorings there is an equality constraint on edges instead of an inequality. This leads to there being no need for the development of skeletons and neighborliness in the strict vector coloring case. However, we do not see how to prove strict vector coloring analogs of Results~\ref{res:Bcor} and~\ref{res:Acor}. This is because we crucially made use of the Perron-Frobenius Theorem for these results, and we cannot apply this to the dual solutions in the strict vector coloring case since we no longer have the nonnegativity requirement. This was one reason why we focused on vector colorings instead of strict vector colorings, since Results~\ref{res:Bcor} and~\ref{res:Acor} are likely to be the useful results in practice.

\paragraph{Acknowledgements:} C.~Godsil is supported by Natural Sciences and Engineering Council of Canada, Grant No.~RGPIN-9439. D.~E.~Roberson is supported by ERC Advanced Grant GRACOL, project no.~320812. R.~\v{S}\'{a}mal is supported by grant 16-19910S of the Czech Science Foundation. A.~Varvitsiotis is supported in part by the Singapore National Research Foundation under NRF RF Award No.~NRF-NRFF2013-13.

\bibliographystyle{plainurl}

\bibliography{UVC3arXiv1.bbl}

\end{document}